\UseRawInputEncoding
\documentclass[11pt]{article}
\usepackage{amsmath, amssymb,ytableau,young,bm,stackrel,textcomp,amsxtra, graphicx}
\usepackage{amsmath, amssymb,amsxtra, graphicx, float,tikz,ytableau,young,bm,stackrel,textcomp,}
\usepackage{bold-extra}
\usepackage{courier}
\usepackage{color}
\usepackage{tabularx}
\usepackage{tikz,pgf}
\usetikzlibrary{shapes.geometric}
\newtheorem{thm}{{\sc Theorem}}[section]
\newtheorem{prop}[thm]{{\sc Proposition}}
\newtheorem{cor}[thm]{{\sc Corollary}}
\newtheorem{lem}[thm]{{\sc Lemma}}
\newtheorem{df}[thm]{{\sc Definition}}
\newenvironment{proof}{\begin{sc}\noindent Proof: \end{sc}}{
     \hbox to 2em{}\nobreak\hfill$\blacksquare$\par\medskip}
\newcommand{\LR}{\hbox{Little\-wood-Richard\-son}}
\newcommand {\PTab}{\hbox{PTab}}

\newcommand {\Pf}{\text{Perf\,}}

\newcommand{\comm}[1]{}
\newcommand {\mmin}{\text{min}}
\newcommand {\mmax}{\text{max}}
\newcommand {\Rot}{\text{Rot}}
\newcommand {\NULL}{\text{NULL}}
\newcommand {\Evac}{\text{Evac}}

\newcommand{\wt}{\text{wt}}
\newcommand{\bw}{\text{Bw}}

\newcommand{\Lus}{\text{Lus}}
\newcommand{\Bw}{\text{Bw}}
\newcommand{\Hw}{\text{Hw}}
\newcommand{\RSK}{\text{RSK}}

\begin{document}

\title{RSK and Ptableaux}
\author{Glenn D.\ Appleby, Tamsen Whitehead\\
Department of Mathematics\\
Santa Clara University\\
Santa Clara, CA 95053\\
{\tt gappleby@scu.edu, tmcginley@scu.edu}}
\date{ }
\newpage

\noindent {\bf Perforated Tableaux in Type $A_{n-1}$ Crystal Graphs and the RSK Correspondence}
\makeatletter
\def\Ddots{\mathinner{\mkern1mu\raise\p@
\vbox{\kern7\p@\hbox{.}}\mkern2mu
\raise4\p@\hbox{.}\mkern2mu\raise7\p@\hbox{.}\mkern1mu}}
\makeatother

\medskip

\noindent
Glenn D. Appleby\\
 Tamsen Whitehead\\
{\em Department of Mathematics\\
and Computer Science,\\
Santa Clara University\\
Santa Clara,  CA 95053}\\
gappleby@scu.edu, tmcginley@scu.edu

\noindent \mbox{} \hrulefill \mbox{}
\begin{abstract} We continue work begun in~\cite{ptab} which introduced \emph{perforated tableaux} as a combinatorial model for crystals of type $A_{n-1}$, emphasizing connections to the classical Robinson-Schensted-Knuth (RSK) correspondence and Lusztig involutions, and, more generally, exploring the role of insertion schemes in the analysis of crystal graphs. An essential feature of our work is the role of \emph{dual} crystals (\cite{GerberLecouvey,vanLeeuwen}) from which we obtain new results within and beyond the classic RSK theory.
\end{abstract}

\section{Introduction}
We regard $\omega \in [n]^{\otimes k}$ as both a word in the alphabet $[n]=\{1,2,\ldots , n\}$ of length $k$, and as a node in the $k$-fold tensor product of the standard crystal of type $A_{n-1}$ (see~\cite{BumpSchilling}). In this paper, we adopt the \emph{Kashiwara convention} to define crystal operators on words so that the map $\omega \mapsto f_i \omega$ puts a crystal structure on $\omega \in [n]^{\otimes k}$ (see below for precise definitions). We extend this to a crystal structure on \emph{biwords} ${\tau \choose \omega}=\binom{\tau_1\tau_2\dots \tau_k}{\omega_1\omega_2\dots\omega_k}$, for $\omega \in [n]^{\otimes k}$ and $\tau \in [m]^{\otimes k}$ (see~\cite{BumpSchilling}, Chapter 9, and \cite{vanLeeuwen}), where we assume if $\tau_{j} = \tau_{j+1}$ then $\omega_{j} \geq \omega_{j+1}$, and set $e_i {\tau \choose \omega} = {\tau \choose e_i \omega}$, producing a crystal structure on biwords that is isomorphic to that on  $[n]^{\otimes k}$. We use Robinson-Schensted-Knuty (RSK) \emph{column} insertion of a biword $\tau \choose \omega$, starting on the \emph{left}, to obtain $(P,Q)$, with $P$ and $Q$ semistandard Young tableaux (SSYT) of the same shape.

Suppose two biwords $\tau \choose \omega$ and $\tau' \choose \omega'$ produce, by RSK insertion, pairs of SSYT $(P,Q)$ and $(P',Q')$, respectively. It is known (see~\cite{BumpSchilling}) that if $P=P'$, then the biwords $\tau \choose \omega$ and $\tau' \choose \omega'$ are \emph{plactically equivalent}, meaning there is an isomorphism $\phi$ of irreducible crystal graphs of biwords such that $\phi{\tau \choose \omega} = {\tau' \choose \omega'}$ (that is, the biwords are in the ``same location'' in the graph). If $Q=Q'$, then the biwords lie in the \emph{same} irreducible crystal graph. Thus, in the context of crystals, classical RSK theory yields a SSYT $P$ as a proxy for the location of a node in an isomorphism class of crystals, and the growth tableau $Q$ becomes a proxy for the specific crystal graph in which a biword is found.

In the authors'~\cite{ptab}, we presented the perforated tableaux (ptableaux) model for type $A_{n-1}$ crystal graphs. Ptableaux are rectangular arrays in which a given box contains a positive integer, or is left blank (``perforated"), with other constraints (defined in the next section).  We defined a bijection ``Perf" from biwords ${\tau \choose \omega}$ to ptableaux $T$ by putting the number $\tau_i$ into row $\omega_i$ of $T$, for $1\le i \le k$. Content in $T$ is arranged in weakly increasing order in rows, and in \emph{strictly} increasing order in columns, by including blank boxes (perforations) as needed to yield column-strictness. We defined crystal operators for ptableaux that were simpler than the standard definitions on words resulting in crystal graphs isomorphic to biword models.

\ytableausetup{centertableaux}

Here, we present a \emph{new} RSK insertion defined on ptableaux, which places classic RSK results in a broader context, and establish new results within the classic theory.  We show ptableau RSK insertion is, in a sense, \emph{dual} to the classic RSK insertion on biwords, but demonstrate several advantages of the ptableau approach.  Among our new results within the classic RSK theory we prove:

\begin{itemize}
\item Not only does the growth tableau $Q$ determine the crystal containing $\tau \choose \omega$, but regarding $Q$ as a perforated tableau we can compute, directly from $Q$, the highest weight element of the irreducible crystal graph containing $\tau \choose \omega$. For example, if
\ytableausetup{smalltableaux}
\[ {\tau \choose \omega} = {1122333444 \choose 2122331331}, \]
then column insertion from the left gives:
\[ \RSK{\tau \choose \omega} = (P,Q)=\left( \,\begin{ytableau} 1&1&1&2&2\\2&3&3&3\\2 \end{ytableau}, \, \begin{ytableau} 1&1&2&3&4\\ 2&3&4&4\\3 \end{ytableau} \, \right). \]
The highest weight element of $\omega = 2122331331$ is determined from $Q$ by writing, from bottom to top, the rows in which the 1's appear in $Q$, then the rows in which the 2's appear, etc. So, in this case, the highest weight of $\omega$ is the word $\eta = 1121321221$, and hence the highest weight of the irreducible crystal containing $\tau \choose \omega$ is $\tau \choose \eta$.
\item Not only does $P$ determine by proxy the location of $\tau \choose \omega$ in its irreducible crystal, but we show how to use $P$ (or, actually, its perforated tableau dual) to compute directly a sequence of crystal operators $e_i$ taking $\tau \choose \omega$ to its highest weight (any two plactically equivalent biwords will have the same path to their respective highest weights). See Theorem~\ref{e up formula} below.
\item Let $\Hw{\tau \choose \omega}$ denote the highest weight element in the irreducible crystal graph containing $\tau \choose \omega$, and let ${\tau \choose \omega} \otimes {\tau_{k+1} \choose \omega_{k+1}}$ denote an \emph{extension} of $\tau \choose \omega$, meaning a valid biword such that:
    \[ \hbox{if }  {\tau \choose \omega} = {\tau_1, \tau_2, \ldots \tau_k \choose \omega_1, \omega_2, \ldots \omega_k},  \hbox{ then } {\tau \choose \omega} \otimes {\tau_{k+1} \choose \omega_{k+1}}= {\tau_1, \tau_2, \ldots \tau_k, \tau_{k+1} \choose \omega_1, \omega_2, \ldots \omega_k, \omega_{k+1}}. \]
We prove: If ${\tau \choose \omega} \otimes {\tau_{k+1} \choose \omega_{k+1}}$ extends $\tau \choose \omega$, then the highest weight $\Hw \left(
{\tau \choose \omega} \otimes {\tau_{k+1} \choose \omega_{k+1}} \right)$ extends the highest weight $\Hw{\tau \choose \omega}$. That is, $\Hw
\left({\tau \choose \omega} \otimes {\tau_{k+1} \choose \omega_{k+1}} \right) = \Hw{\tau \choose \omega} \otimes {\tau_{k+1} \choose \omega_{k+1}'}$ for some $\omega_{k+1}' \in [n]$.
\end{itemize}

We prove these results by re-interpreting the classic RSK theory as one that is \emph{dual} (in a sense defined below, see~\cite{ptab}) to our RSK insertion on perforated tableaux. We define an insertion scheme on a ptableau $T$, producing a pair of ptableaux $(PT, T_{\mmax})$, and prove:

\begin{itemize}
\item If $\RSK(T) = (PT,T_{\mmax})$, then $T_{\mmax}$ is the highest weight element in the irreducible crystal containing $T$. That is, the insertion computes the highest weight element \emph{directly}, and not by proxy.
\item The ptableau $PT$ is an element of a uniquely determined \emph{distinguished crystal}, characterized by the property that its highest weight element is a SSYT that has only $i$'s in row $i$.
\item If $\RSK(T) = (PT, T_{\mmax})$, then $T$ is \emph{plactically equivalent} to $PT$. In particular, if $\RSK(T) = (PT, T_{\mmax})$ and $\RSK(T') = (PT, T_{\mmax}')$, then $T$ and $T'$ are plactically equivalent (since they are equivalent to a common ptableau $PT$ in the distinguished crystal). That is, we construct explicit plactic equivalences (to a common element in the distinguished crystal).
\item Suppose $\Pf{\tau \choose \omega} = T$ under the bijection between biwords and ptableaux mentioned above. Then, if $\RSK{\tau \choose \omega} = (P,Q)$ and $\RSK(T) = (PT, T_{\mmax})$, we prove $T_{\mmax} = Q$. That is, under the \emph{classic} column insertion (starting on the left) on biwords, the growth tableau $Q$ actually \emph{is} the highest weight element of the perforated tableau determined by the biword $\tau \choose \omega$.
\item Finally, it was proved in~\cite{ptab} that SSYT, when viewed as ptableaux, are \emph{highest weight elements}. If $\RSK{\tau \choose \omega} = (P,Q)$ and $\RSK(T) = (PT, T_{\mmax})$, then regarding $P$ as a \emph{ptableau}, we prove $PT = \widehat{P}$, where $\widehat{P}$ is the \emph{dual} of $P$ (see Section 2.3). In particular, elements of the distinguished crystal (such as $PT$) are dual to highest weight ptableaux.
\end{itemize}

In short, we show that results of classic RSK theory on biwords, in the context of crystal graphs, can be viewed within the crystal theory of perforated tableaux, where we achieve our results more simply and directly. The classic RSK map on biwords $RSK{\tau \choose \omega} = (P,Q)$ yields SSYT $P$ and $Q$ that are distinct from the crystal model on biwords. We feel it is more natural to use the map $RSK(T)=(PT,T_{\mmax})$ in which all objects are ptableaux and from which we can compute information regarding this invariant directly from the RSK output, and not by proxy.

\section{Definitions, Perforated Tableaux}

Definitions and basic results regarding perforated tableaux can be found in \cite{ptab}. We give a brief overview here, along with a short review of crystal operators on words.

\subsection {Crystal Operators on Words and Biwords}

Let $[n] = \{1,2, \ldots , n\}$ and let $[n]^{\otimes k}$ denote the length-$k$ words in the alphabet $[n]$. We regard $[n]^{\otimes k}$ as the $k$-fold tensor product of the standard crystal graph $[n]$ of type $A_{n-1}$ (see~\cite{BumpSchilling}). In this paper we adopt the \emph{Kashiwara Convention} for crystal operators on $[n]^{\otimes k}$:

\begin{df}Suppose $\omega = \omega_{1}\omega_{2} \cdots \omega_k \in [n]^{\otimes k}$. For any $s \in [n]$, define $\phi_{i}(s) = \delta_{i}(s)$ (Kronecker delta function), and $\epsilon_{i}(s) = \delta_{i+1}(s)$. Then let

\begin{align*} c_{e}^{i}(\omega,j)  &=   \left( \sum_{s=1}^{j}\epsilon_{i}(\omega_{s}) - \sum_{s =1}^{j-1}\phi_{i}(\omega_{s}) \right),\\
c_{f}^{i}(\omega,j)&=\left( \sum_{s=j}^{k} \phi_{i}(\omega_{s}) - \sum_{s = j+1}^{k} \epsilon_{i}(\omega_{s}) \right), \end{align*}
and, for $1 \leq i <n$,
\[ \epsilon_{i}(\omega) = \max_{j=1}^{k}c_{e}^{i}(\omega,j), \qquad \phi_{i}(\omega) = \max_{j=1}^{k} c_{f}^{i}(\omega,j). \]

Then for $s\in [n]$, define $e_{i}(s) = i$ if $s= i+1$, and $\NULL$ otherwise, and $f_{i}(s) = i+1$ if $s=i$, and $\NULL$ otherwise. Extend these definitions to the word $\omega \in [n]^{\otimes k}$ by setting
\[ e_{i}(\omega) = \omega_{1} \omega_{2}  \cdots \omega_{j-1} e_{i}(\omega_{j})\omega_{j+1} \cdots \omega_k, \]
with $j$ the \emph{smallest} index at which $\epsilon_{i}(\omega) = c_{e}^{i}(\omega,j)$, and
\[ f_{i}(\omega) = \omega_{1} \omega_{2} \cdots \omega_{j-1}  f_{i}(\omega_{j})\omega_{j+1} \cdots \omega_k, \]
with $j$ the \emph{largest} index at which $\phi_{i}(\omega) =  c_{f}^{i}(\omega,j)$, or $\NULL$ if $e_{i}(\omega_{j})$ or $f_{i}(\omega_{j})$ is $\NULL$.

The \emph{weight} of a word $\omega \in [n]^{\otimes k}$ is the composition $\wt(\omega) = (a_1, \ldots , a_n)$ where $a_i$ is the number of $i$'s in $\omega$.
\end{df}

\begin{prop}
The definitions given above determine an $A_{n-1}$ crystal structure on $[n]^{\otimes k}$, as a $k$-fold tensor product of the standard crystal $[n]$. \label{crystal def}
\end{prop}

As stated in the previous section, we can extend the crystal structure on words to \emph{biwords}.

\begin{df} (1) Suppose $\omega \in [n]^{\otimes k}$ and $\tau \in [m]^{\otimes k}$. Let $\omega \choose \tau$ denote a \emph{biword}:
\[ {\tau \choose \omega} = {\tau_1, \tau_2, \ldots \tau_k \choose \omega_1, \omega_2, \ldots \omega_k}. \]

(2) A biword $\tau \choose \omega$ is in \emph{standard form} if the letters of $\tau$ are weakly increasing, and the letters of $\omega$ under any constant factor of $\tau$ are weakly decreasing. That is, if $\tau_i = \tau_{i+1}$ then $\omega_{i} \geq \omega_{i+1}$.

(3) Let $\Bw([m]\times [n])^{\otimes k}$ denote the set of biwords of length $k$ in the alphabets $[m]$ and $[n]$. That is, given any subset $S$ (of size $k$) of the cartesian product $[m]\times [n]$ (regarded as a set of ordered pairs of the form $ \tau_i\choose \omega_i$, with $\tau_i \in [m]$ and $\omega_i\in [n]$), then ${\tau \choose \omega} \in \Bw([m] \times [n])^{\otimes k}$ is the uniquely determined biword obtained by ordering the elements of $S$ so that $\tau \choose \omega$ is in standard form.
\end{df}

We put a type $A_{n-1}$ crystal structure on biwords by defining
\[ e_i {\tau \choose \omega} =  {\tau \choose e_i\omega} \]
using the action of $e_i$ on $\omega \in [n]^{\otimes k}$ given above, and setting the weight of $\tau \choose \omega$ to be the weight of the bottom row $\omega$. In~\cite{ptab}, Theorem 5.1, it was shown that this action preserves the biword structure in that if $\tau \choose \omega$ is in standard form, then so is ${\tau \choose e_i\omega}$.

\subsection {Perforated Tableaux}
In~\cite{ptab}, the combinatorial model of \emph{perforated tableaux} was introduced. These are objects such as:

\[ T=\begin{ytableau} \ & & 1 && 3 & 4\\ & 1 &2 & 2&& \\3 & 3 & 4 & 4 && \end{ytableau} \]
that form the nodes of type $A_{n-1}$ crystal graphs. More formally:

\begin{df} \label{ptab} A rectangular tableau $T$ is called a \emph{perforated tableau}  (or \emph{ptableau}) if:
 \begin{enumerate}
 \item All boxes in $T$ contain positive integers (have content), or are blank (unfilled, denoted $\Box$).
  \item For any positive integer $i$, the $i$'s in $T$ form a horizontal strip.
 \item If $i < j$, then no element in the horizontal strip of $j$'s is in a higher row, nor in a column to the left, of a column containing any element of the horizontal strip of $i$'s.
 \item The number of columns of $T$ subject to the above conditions is \emph{minimum}. In particular, $T$ has no column entirely composed of blanks (rows of all blanks are allowed).
\end{enumerate}

If a box in a ptableau $T$ is not a blank, we say it has \emph{content}, and refer to the boxes of $T$ with content as the \emph{cells} of $T$.

Given a perforated tableau $T$, we define the \emph{weight} of $T$ as the composition $\wt(T) = (a_1, \ldots , a_n)$ where $a_i$ is the number of cells in row $i$ of $T$.
\end{df}

\noindent \emph{Notational Convention:} We need to distinguish between an actual box of a ptableau, and the \emph{value} (an integer or blank) of that box.  We denote a cell with content $a$ by $\framebox{a}$ and a blank box by $\Box$. We denote a particular $\framebox {a}$ in row $t$ by $\framebox{a}_{\, t}$.

A given perforated tableau can have several equivalent forms:

\begin{df} Suppose $S$ and $T$ are ptableaux with $n$ rows. Then $S$ and $T$ are \emph{row-equivalent} if we can transform $S$ into $T$ by a finite number of steps in which a blank in a given row is swapped with an adjacent cell (non-blank) in the same row, such that the result is a valid ptableau. Note that if $S$ and $T$ are row-equivalent, then $\wt(S) = \wt(T)$.
\end{df}
Two particular row-equivalent ptableaux will be used extensively in subsequent sections:
\begin{df}
A ptableau $T$ is \emph{left-justified} if all of $T$'s content is as far to the \emph{left} as possible. That is, in any ptableau $T'$ that is row-equivalent to $T$, each entry of $T$ is weakly to the left of the corresponding entry in $T'$.

Similarly, a ptableau $T$ is \emph{right-justified} if all content of $T$ is as far right as possible.
\label{justified}
\end{df}

The following are then easily proved:

\begin{lem}
\begin{enumerate}
\item Given any ptableau $T$, there is a unique ptableau $^*T$ that is left-justified and row-equivalent to $T$.
\item Given a left-justified ptableau $^*T$, there is a unique right-justified ptableau $T^*$ that is row-equivalent to $^*T$.
\end{enumerate}\label{left-right-justified}
\end{lem}

So, for example, if

\[ T= \begin{ytableau}
\ &&&&1&1&4 \\
&1&&1&&2&5\\
1&2&3&4&4&5& \end{ytableau}\, , \]
then we have the row-equivalent ptableaux:
\[ ^*T=\begin{ytableau}
\ &&&1&1&4& \\
&1&1&2&&&5\\
1&2&3&4&4&5& \end{ytableau} \ \ \hbox{and} \   T^* = \begin{ytableau}
\ &&&&1&1&4 \\
&&1&1&&2&5\\
1&2&3&4&4&5&\end{ytableau}\, . \]

A consequence of Lemma 3.6 of ~\cite{ptab} is that all perforated tableaux in a given row-equivalence class have the same number of columns. Row-equivalence is also clearly an equivalence relation on ptableaux.

\begin{df} Let $\PTab$ denote the set of row-equivalence classes of ptableaux, each of which contains ptableaux that are row-equivalent to the left-justified representative $^*T$ and that have the same number of columns as $^*T$.  Let $\PTab_{n}$ denote the subset of $\PTab$ composed of ptableaux with $n$ rows (some of those rows may be entirely blank).

Let $\PTab_{(m,n)}$ be the subset of $\PTab_n$ with content from $[m]$, and $\PTab_{(m,n)}^k$ denote the subset of $\PTab_{(m,n)}$ with $k$ cells (filled boxes).
\label{ptab def}
\end{df}

Again, the following is easily proved:

\begin{lem} A given ptableau $T \in \PTab$ is determined by its content, and the rows in which that content appears. \label{unique}
\end{lem}
With this, we have
\begin{df} The map
\[ \Pf: \Bw([m]\otimes [n])^{\otimes k} \rightarrow \PTab_{(m,n)}^k \]
is defined by setting $\Pf {\tau \choose \omega}$ to be the ptableau uniquely determined from $\tau \choose \omega$ (by Lemma~\ref{unique}) with content $\tau_i$ in row $\omega_i$, for each $i$, $1 \leq i \leq k$.
\end{df}

It was proved in ~\cite{ptab}, Lemma 3.5, that $\Pf: \Bw([m]\otimes [n])^{\otimes k} \rightarrow \PTab_{(m,n)}^k$ is a weight-preserving bijection. We denote the inverse map by $\Bw: \PTab_{(m,n)}^k \rightarrow  \Bw([m]\otimes [n])^{\otimes k}$ with $\bw(T)$ the biword containing a column $\tau_i \choose \omega_i$ for each $\tau_i$ in row $\omega_i$ of $T$, for $1 \leq i \leq k$, written in standard form.

We now define a crystal structure on ptableaux.

\begin{df}
Let $T \in \PTab_{n}$  and choose $i$, for $1 \leq i <n$. Let $T[i,i+1]$ denote the ptableau formed by the entries of $T$ in rows $i$ and $i+1$.\end{df}

Note that, as columns of blanks are not allowed in a ptableau,  we omit any such columns in the construction of $T[i,i+1]$. For example, if
\ytableausetup{smalltableaux}
\[ T=  \begin{ytableau}
\ &&&&&&&&4&4&5\\
&&&&&1&1&2&&&6\\
&&&1&1&&2&4&5&6&7\\
1&1&2&3&3&3&4&6&6&& \end{ytableau} \, , \quad \hbox{then} \ \  T[2,3] = \begin{ytableau}
\ &&1&1&2&&6\\
1&1&2&4&5&6&7 \end{ytableau} \, . \]

Lemma 4.6 in ~\cite{ptab} shows that if content in row $(i+1)$ of $^*T[i,i+1]$ has a blank above it, then the corresponding content in $^*T$ also has a blank above it. With this, we have
\begin{df} \label{e_i def} Suppose $T \in \PTab_{n}$ and that $\framebox{c}_{\, i+1}$ is an entry in row $i+1$ of $T$ such that in $^*T[i,i+1]$ (the left-justified form of $T[i,i+1]$),  the corresponding $\framebox{c}$ is the right-most entry with a blank above it. Then $e_i (T)$ is the ptableau obtained by swapping $\framebox{c}_{\, i+1}$ with the blank above it (in row $i$) of $T$. If there is \emph{no} such entry $\framebox{c}_{\, i+1}$, then $e_i (T) = \NULL$.

\medskip
\noindent We define $\epsilon_{i}(T)$ to be the number of blanks in row $i$ (the top row) of $^*T[i,i+1]$.
\end{df}

\begin{df}  \label {f_i def} Suppose $T \in \PTab_{n}$ and that $\framebox{c}_{\, i}$ is an entry in row $i$ of $T$ such that in $T[i,i+1]^*$ (the right-justified form of $T[i,i+1]$) the corresponding $\framebox{c}$ is the left-most entry with a blank below it. Then,  define $f_i (T)$ to be the ptableau obtained by swapping $\framebox{c}_{\, i}$ with the blank below it (in row $i+1$) of $T$. If there is \emph{no} such entry $\framebox{c}$, then $f_i (T) = \NULL$.

\medskip

\noindent We define $\phi_{i}(T)$ to be the number of blanks in row $i+1$ (the bottom row) of $T[i,i+1]^*$.
\end{df}

In~\cite{ptab} it was shown the operations $e_i$ and $f_j$ above put a type $A_{n-1}$ crystal structure on $\PTab_{n}$. As an example, consider $T$ and $T[2,3]$ as above.

We calculate $e_{2}(T)$:
\begin{equation*} ^*T[2,3] = \begin{ytableau}
\ &&1&1&2&&6\\
1&1&2&4&5&6&7 \end{ytableau} ,\label{example eq} \end{equation*}
and so the $6$ in the sixth column of $^*T[2,3]$ is the rightmost entry under a blank. We circle this entry in $^*T$ below:
\ytableausetup{centertableaux}
\[ ^*T=  \begin{ytableau}
\ &&&&&&&&4&4&5\\
&&&&1&1&2&&&&6\\
&&1&1&2&&&4&5&\textcircled{$\scriptstyle 6$}&7\\
1&1&2&3&3&3&4&6&6&& \end{ytableau}\, ,
 \hbox{ and so }
e_{2}(T) =  \begin{ytableau}
\ &&&&&&&&4&4&5\\
&&&&1&1&2&&&\textcircled{$\scriptstyle 6$}&6\\
&&1&1&2&&&4&5&&7\\
1&1&2&3&3&3&4&6&6&& \end{ytableau}\, . \]
By definition, $\epsilon_2(T) = 3$ because there are three blanks in the top row of $T[2,3]$.
 As proven in~\cite{ptab} we have:

\begin{thm}[\cite{ptab}]
Under the bijection $\Pf$,  the type $A_{n-1}$ crystal structure on $\PTab_{(m,n)}^k$ is isomorphic to that obtained on biwords  $\Bw([m]\otimes [n])^{\otimes k}$.\label{biword iso}
\end{thm}

Thus, computing crystal operators on ptableaux is especially simple. In particular,  Theorem 6.2 of~\cite{ptab} states:

\begin{thm} A ptableau $T$ is highest weight if and only if its left-justified form $^*T$ is a semistandard Young tableau. \label{highest weight ptab}
\end{thm}

For example, the ptableau
\[ T = \begin{ytableau}
\ &&&&1&1&4 \\
&&1&1&&2&5\\
1&2&3&4&4&5& \end{ytableau} \ \
\hbox{is not highest weight since} \
 ^{*}T =  \begin{ytableau}
\ &&&1&1&4 &\\
&1&1&2&&&5\\
1&2&3&4&4&5& \end{ytableau} \, , \]
which is not a SSYT. However,
the ptableau
\[ T' =  \begin{ytableau}
 1&1&1&1&1&4&4 \\
& 2 && 4&4&6&\\
&&3&&5&& \end{ytableau} \, , \ \
\hbox{is highest weight since}\
^{*}T' = \begin{ytableau}
 1&1&1&1&1&4&4 \\
 2 & 4&4&6&&&\\
3&5&&&&& \end{ytableau}\, , \]
which is semistandard.

\begin{df} We say two biwords ${\tau \choose \omega}$ and ${\tau' \choose \omega'}$ are \emph{plactically equivalent} if there is an isomorphism $\phi$ of irreducible, connected crystal graphs of biwords such that $\phi{\tau \choose \omega} = {\tau' \choose \omega'}$. We denote this by ${\tau \choose \omega}
\equiv {\tau' \choose \omega'}$. Similarly, we say two ptableaux $T$ and $T'$ are plactically equivalent, denoted $T \equiv T'$ if $\bw(T) \equiv \bw(T')$.
\end{df}

We will need the following notational definitions:

\begin{df}
(1) Given a biword $\tau \choose \omega$, let ${\tau \choose \omega}^{(\ell)}$ denote the biword formed by the first $\ell$ many columns of $\tau \choose \omega$. Similarly, let $T^{(\ell)}$ denote the ptableau formed by $\bw(T)^{(\ell)}$.

(2) Given a biword ${\tau \choose \omega} \in \Bw([m]\otimes [n])^{\otimes k}$,
\[ {\tau \choose \omega} = {\tau_1, \tau_2, \ldots \tau_k \choose \omega_1, \omega_2, \ldots \omega_k}, \]
we denote by ${\tau \choose \omega} \otimes {\tau_{k+1} \choose \omega_{k+1}} \in \Bw([m]\otimes [n])^{\otimes k+1}$ the \emph{extension} of $\tau \choose \omega$ by the column $\tau_{k+1} \choose \omega_{k+1}$ if
\[{\tau \choose \omega} \otimes {\tau_{k+1} \choose \omega_{k+1}}= {\tau_1, \tau_2, \ldots \tau_k, \tau_{k+1} \choose \omega_1, \omega_2, \ldots \omega_k, \omega_{k+1}} \]
is in standard form. Similarly, let $T \otimes \framebox{$\tau_{k+1}$}_{\, \omega_{k+1}}$ denote the ptableau $T$, with $\tau_{k+1}$ appended to row $\omega_{k+1}$, assuming this is a valid ptableau.  Then $T \otimes \framebox{$\tau_{k+1}$}_{\, \omega_{k+1}}$ is an \emph {extension} of $T$ provided $\bw(T \otimes \framebox{$\tau_{k+1}$}_{\, \omega_{k+1}}) = {\tau \choose \omega} \otimes {\tau_{k+1} \choose \omega_{k+1}}$ is an extension of $\bw(T) = {\tau \choose \omega}$.
\label{bi-word def}
\end{df}

For example, if
\[ T =\begin{ytableau} \ & & 1 & 3 & 4 \\ & 1 & 2 & 2 & \\ 3 & 3 & 4 & 4& \end{ytableau}, \]

then the biword of $T$ is:
\[\bw(T)={
1122333444\choose
2122331331}
,\]
and
\[ \bw(T)^{(5)} ={
11223\choose
21223} \ \ \hbox{and} \ \  T^{(5)} = \begin{ytableau} \ &  1 &    \\  1 & 2 & 2  \\ 3 &  &    \end{ytableau}. \]

\begin{df} Given a ptableau $T$, we say a subset of cells of $T$ form a \emph{right-contiguous strip} if

\begin{enumerate}
\item The cells form a \emph{right-horizontal strip}, meaning that in $T^*$, the right-justified form of $T$, no cells in the strip appear in the same column, and if $\framebox{$s$}_{\, k}$ is a cell in the strip lying in a column to the left of a cell $\framebox{$s'$}_{\, k'}$ in the strip, then $k \geq k'$.
\item In the right-justified form $T$, if two cells of the strip appear in columns $j$ and $j'$, $j < j'$, then there are cells in the strip appearing in all columns $j''$, for $j \leq j'' \leq j'$.
\end{enumerate}
\end{df}

\begin{df} Define a map $\Rot : \PTab_{(m,n)} \rightarrow \PTab_{(m,n)}$ with $\Rot(T)$ the ptableau obtained from $T$ by first rotating $T$ by $180^\circ$, and then replacing each number $j$ by $m-j+1$.\label{rot def}
\end{df}

For example, if
\ytableausetup{smalltableaux}
 $T =  \begin{ytableau}
\ &&&&&&&&4&4&5\\
&&&&&1&1&2&&5&6\\
&&&1&1&&2&4&5&6&7\\
1&1&2&3&3&3&4&6&6&&8 \end{ytableau}, $
then
$ \Rot(T) = \begin{ytableau}
\ 1 & & 3 &3&5&6&6&6&7&8&8 \\
2&3&4&5&7&&8&8&&&\\
3&4&&7&8& 8&&&&&\\
4&5&5&&&&&&&& \end{ytableau}. $

\medskip

The following in easy to show:
\begin{prop}
Given $T \in \PTab_{(m,n)}$, we have

\[ e_i\Rot(T) = \Rot(f_{n-i}T). \]

In particular, if $T$ and $T'$ are plactically equivalent, then so are $\Rot(T)$ and $\Rot(T')$. \label{rot prop}
\end{prop}

\bigskip

\subsection{Duality}

We now define the \emph{dual} of a biword and ptableau (see \cite{ptab}, Section 1.2), and state some associated results.

\begin{df} Given a biword $\tau \choose \omega$, the \emph{dual} of $\tau \choose \omega$, denoted $\widehat{\tau \choose \omega}$, is the biword ${\widehat{\omega} \choose \widehat{\tau}}$ formed by putting the biword $\omega \choose \tau$ into standard form.
\end{df}

This definition allows us to put a $GL_{n} \times GL_{m}$ \emph{bicrystal structure} on biwords  $\Bw([m]\otimes [n])^{\otimes k}$. Let $e_{i}^{(n)}$ denote an operator for the $GL_n$ structure on words $\omega \in [n]^{\otimes k}$, and similarly let $e_{j}^{(m)}$ denote an operator for the $GL_m$ structure on words $\tau \in [m]^{\otimes k}$. Then the bicrystal structure on biwords $\Bw([m]\otimes [n])^{\otimes k}$ (see~\cite{BumpSchilling}, Chapter 9)  is given by setting:
\[ e_{i}^{(n)} {\tau \choose \omega} = {\tau \choose e_{i}^{(n)} \omega} \ \ \hbox{and} \ \ e_{j}^{(m)}{\tau \choose \omega} = \widehat{\widehat{\omega} \choose e_{j}^{(m)} \widehat{\tau}}. \]

We review some facts regarding duality and RSK insertion. Recall that, in this paper, by ``RSK insertion" of a biword $\tau \choose \omega$ we mean column insertion of $\omega$, starting on the \emph{left} of the biword, with growth recorded by $\tau$. Results stated here would, in some form, remain true under other insertion schemes, but the form we adopt is most easily related to ptableaux.

Let $RSK{\tau \choose \omega} = (P,Q)$ be the result of RSK insertion of the biword $\tau \choose \omega$, producing a SSYT pair $(P,Q)$ of the same shape. Kashiwara's tableau model~\cite{Kash-Nak} defines a crystal operator on SSYT: $P \mapsto e_i P$ and a crystal isomorphism from SSYT pairs to biwords such that if $RSK{\tau \choose \omega} = (P,Q)$, then ${\tau \choose \omega} \mapsto {\tau \choose e_i \omega}$ implies $P \mapsto e_i P$. It is also classically known that $RSK{\tau \choose \omega} = (P,Q)$, then $RSK\widehat{\tau \choose \omega} = (Q,P)$.

It has been shown (see~\cite{BumpSchilling}) that two biwords $\tau \choose \omega$ and $\tau' \choose \omega'$ are plactically equivalent if and only if they have the same $P$ tableau under RSK insertion. Similarly, they lie in the same irreducible crystal graph when they have the same growth tableau $Q$. From these observations we see that two biwords lie in the same irreducible crystal graph if an only if their duals are plactically equivalent, summarized by the commutative diagram:

\[ \begin{array}{cccl} (P,Q) & \stackrel{e_i}{\rightarrow} & (e_{i}P,Q)& \hbox{$\Leftarrow$ In the same crystal} \\ \raisebox{-0.05in}{$\widehat{\ }$} \downarrow &&\raisebox{-0.05in}{$\widehat{\ }$} \downarrow& \\
(Q,P) & \stackrel{e_i}{\rightarrow} & (Q,e_{i}P) & \hbox{$\Leftarrow$ Plactically equivalent} \end{array} \]

\noindent We construct ptableaux analogues for these observations.

\begin{df} Given a ptableau $T \in \PTab_{(m,n)}^k$ we form the \emph{dual} ptableau $\widehat{T} \in \PTab_{(n,m)}^k$  by requiring that $\widehat{T}$ has a cell with content $j$ in row $i$ for every cell in $T$ with content $i$ in row $j$.
\label{dual def}
\end{df}

\[ \hbox{So, if} \ \ T =\begin{ytableau} \ & & 1 & 3 & 4 \\ & 1 & 2 & 2 & \\ 3 & 3 & 4 & 4& \end{ytableau} \, , \ \ \hbox{then} \ \ \widehat T=\begin{ytableau} \ & 1&2 & 2\\ 1 & &3 & 3\\ 2&3&& \end{ytableau} \, . \]

With this definition, since it essentially reverses the role of row index and content, it is easy to see that duality (for ptableaux and biwords) commutes with the map $\Pf$:
\[ \Pf\, \widehat{\tau \choose \omega} = \widehat{\Pf {\tau \choose \omega}}. \]

We can produce $\widehat{T}$, the dual of $T$, directly from the biword $\bw(T)={\tau \choose \omega}$: put the factor of $\omega$ that is under the $1$'s of $\tau$, in reverse order, into row $1$ of $\widehat{T}$, the factor  of $\omega$ that is under the $2$'s of $\tau$, in reverse order, in row $2$, etc., inserting blanks as needed to ensure column strictness. For example,
\[ \hbox{If, as before,} \ \ T =\begin{ytableau} \ & & 1 & 3 & 4 \\ & 1 & 2 & 2 & \\ 3 & 3 & 4 & 4& \end{ytableau}, \ \ \hbox{then} \ \ \bw(T)={\tau \choose \omega} = { 11122233 \choose 22133132}. \]
The top row of 1's, 2's and 3's in $\tau$ breaks $\omega$ into the factors $221$, $331$, and $32$. To compute $\widehat T$ from this, we reverse these three factors of $\omega$ and put them in rows 1, 2 and 3, with blanks as needed for column-strictness:
\[ \widehat T=\begin{ytableau} \ & 1&2 & 2\\ 1 & &3 & 3\\ 2&3&& \end{ytableau}. \]

The simplicity of the map from biwords to dual ptableaux might suggest that duals $\widehat T$ would be a better crystal model than the images $\Pf{\tau \choose \omega} = T$ for crystal graphs. However, ptableaux have several features that make the analysis of crystal graphs far easier, as evidenced, for example, in Theorem~\ref{highest weight ptab} above.

In addition to biwords, and ptableaux, there are \emph{matrix models} for crystal graphs~\cite{BumpSchilling,Fulton,Shimozono,vanLeeuwen}. We compute a  \emph{matrix} $M$ associated to a biword $\tau \choose \omega$, by setting $\alpha_{ij}$, the $(i,j)$ entry of $M$, equal to the number of times $i$ appears over $j$ in the biword. This map is a bijection from biwords ${\tau \choose \omega} \in \Bw([m]\otimes [n])^{\otimes k}$ to $\text{Mat}_{m \times n}({\mathbb N})^{(k)}$, the $m \times n$ matrices with non-negative integer entries with the property that the sum of the matrix entries is $k$. In this case, a $GL_n$ crystal structure can be defined for matrices $M \in \text{Mat}_{m \times n}({\mathbb N})^{(k)}$ by computing along the \emph{columns} of $M$, and one shows this structure commutes with the map from biwords to matrices, yielding isomorphic crystal structures. We obtain a bicrystal structure by defining the corresponding $GL_{m}$ structure on the \emph{transpose} of $M$, so that the dual of $M$ is $\widehat{M} = M^{T}$.

\medskip

Duality is easily analyzed in the ptableaux setting. In fact, we can determine much of the structure of $\widehat{T}$ directly from $T$. An example of this, and a necessary result for our work below, is the following:

\begin{thm} Let $T \in \PTab_n$, and suppose $T$ has $\ell$-many columns. Then column $\ell-k+1$ of $T^*$ (the right-justified form of $T$) is dual to column $k$ of the $^*\widehat T$ (the left-justified form of the dual $\widehat{T}$). \label{left-right}
\end{thm}

\begin{proof}
Recall that the notation $\framebox{s}_{\, t}$ denotes a particular cell $\framebox{s}$ in row $t$ of a ptableau.  We call the horizontal strip of $i$'s in a ptableau $T$ the \emph{$i$-strip} of $T$. We begin with column $\ell$, the right-most column of $T^*$. If a cell $\framebox{$s$}_{\, t}$ is in this column, then this cell is the end of the $s$-strip of $T^*$, so that $t$ is minimal among all such row indices with content $s$. Hence, by definition of the dual of $T$, there is a $t$ in row $s$ of $\widehat{T}$ , which we call $\framebox{$t$}_{\, s}$. Since the content, and their corresponding row indices,  in the right-most column of $T$ are strictly increasing, the cell $\framebox{$t$}_{\, s}$ appears in the left-most position, and hence is the \emph{start} of the $t$-strip in $\widehat T$, so that it must appear in the left-most column of $^*\widehat{T}$.

We assume, inductively, that the conclusions of the theorem hold in columns $\ell -j+1$, for $1 \leq j < k$. We prove they hold in column $\ell-k+1$ as well. If some cell $\framebox{$s_k$}_{\, t_k}$ appears in column $\ell-k+1$ of $T^*$, then there must be a right-contiguous strip of cells
$\framebox{$s_{j}$}_{\, t_{j}}$, $1 \leq j \leq k$ in columns $n-j+1$ such that
$s_j \geq s_{j+1}$, and $t_j \leq t_{j+1}$.

That is, as $j$ decreases from $k$ to $1$ (moving left to right in $T$), there is a right-contiguous strip of cells with content weakly increasing, rows weakly decreasing, that prevents us from moving the cell $\framebox{$s_k$}_{\, t_k}$ any farther to the right.

Dually, this implies that there are cells  $\framebox{$t_{j}$}_{\, s_{j}}$, $1 \leq j \leq k$, in columns $j$ of the left-justified form of $\widehat{T}$, such that the content is weakly increasing (left to right), with the rows weakly decreasing, and hence forming a left-contiguous strip in the left-justified form of $\widehat{T}$, so that the cell $\framebox{$t_{k}$}_{\, s_{k}}$ of must appear in column $k$, as claimed.
\end{proof}

Thus, for example, we have:

\[ T =\begin{ytableau} \ &&1&3\\ 1&1&&\\2&2&3& \end{ytableau}, \quad  \hbox{and} \quad  \widehat T=\begin{ytableau} \ & 1&2 & 2\\  & &3 & 3\\ 1&3&& \end{ytableau}. \]

We see that in column $k=2$ of the right-justifed $T$ we have a $1$ in row $2$ and a $2$ in row $3$ so that, dually, in column $3$ of the left-justified $\widehat T$ we therefore have a $2$ in row $1$ and a $3$ in row $2$.

Note that, typically, a ptableau $T$ and its dual $\widehat T$ need not have the same number of rows, though by the above result, they must have the same number of \emph{columns}.

\medskip

\noindent The following will be necessary in our ptableau-inspired version of RSK:

\begin{df} Let $T^*$ be the right-justified form of a ptableau $T$. Then a cell $\framebox{$s$}_{\, k}$ in $T^*$ is a \emph{violation entry} if there is no cell with content $s-1$ in the same column of $T^*$ as $\framebox{$s$}_{\, k}$. If the first cell with content above $\framebox{$s$}_{\, k}$ has some value $s'$, with $s' < s-1$, then the violation entry $\framebox{$s$}_{\, k}$ has \emph{multiplicity} $s - s' -1$ (if there is \emph{no} entry above the content $s$ set $s'=0$). \label{violation entry}
\end{df}

\begin{cor} If a ptableau $T$ has a violation entry $\framebox{$s$}_{\, k}$ in column $\ell - j+1$ of multiplicity $m$ in its right-justified form $T^*$, then the left-justified form of the dual $\widehat T$ has an entry $\framebox{$k$}_{\, s}$ with $m$ blanks above it in column $j$.\label{blank dual cor}
\end{cor}
\begin{proof} If $\framebox{$s$}_{\, k}$ is a violation entry of multiplicity $m$ in column $\ell-j+1$, then there are \emph{no} cells above $\framebox{$s$}_{\, k}$ in this column with values $s -1, s-2, \ldots , s-m$. Hence, by Theorem~\ref{left-right} there are no cells in column $j$ of $^*\widehat{T}$ (the left-justified form of $\widehat T$) appearing in rows $s-1, \ldots , s-m$, proving our claim.
\end{proof}

For example, the circled entries in $T^*$ below are violation entries;  the $3$ in row $1$ has multiplicity $3-0-1=2$ and the $3$ in row $3$ has multiplicity $3-1-1=1$.

\[ T^* =\begin{ytableau} \ &&1&\textcircled{$\scriptstyle \bf 3$}\\ 1&1&&\\2&2&\textcircled{$\scriptstyle \bf 3$}& \end{ytableau}.\]
We circle the corresponding entries in $\widehat T$:

\[ ^*\widehat T=\begin{ytableau} \ & 1&2 & 2\\  & &3 & 3\\ \textcircled{$\scriptstyle \bf 1$}&\textcircled{$\scriptstyle \bf 3$}&& \end{ytableau}. \]
The two blanks over the \textcircled{$\scriptstyle \bf 1$} in row $3$ of $^*\widehat T$ correspond to the violation entry \textcircled{$\scriptstyle \bf 3$} of multiplicity $2$ in row $1$ of $T^*$.  Likewise, the blank over the \textcircled{$\scriptstyle \bf 3$} in row $3$ of $^*\widehat T$ corresponds to the violation entry \textcircled{$\scriptstyle \bf 3$} of multiplicity $1$ in row $3$ of $T^*$.

\subsection{The Distinguished Crystal}

\begin{df}
Given some partition $\mu=(\mu_1, \mu_2,\dots, \mu_n)$, let $T_\mu$ be the ptableau of shape $\mu$, with $\mu_1$ $1$s in row $1$, $\mu_2$ $2$s in row $2$, etc.  Then the irreducible crystal ${\cal B}\in \PTab_n$ whose highest weight is $T_\mu$ is called the \emph {distinguished crystal} of highest weight $\mu$.\label{distinguished crystal def}
\end{df}

\begin{df} A ptableau $T \in \PTab_n$ is said to satisfy the \emph{word condition} if for all $i$, $1 \leq i \leq n-1$, the number of $i$'s in rows $i$ through $i+k$ of $T$ is greater than or equal to the number of $(i+1)$'s in rows $i+1$ through $i+k+1$.
\end{df}

Note that the word condition on ptableaux does \emph{not} correspond to the Yamanouchi condition on words $\omega$, even if $\Pf(\omega) = T$. The name is chosen only to be consistent with usage of the term in characterizing \LR\ tableaux.

\begin{cor} The following are equivalent for a ptableau $T$

 \begin{enumerate}
\item $T$ satisfies the word condition.
\item $T$ lies in the distinguished crystal of highest weight $\mu = (\mu_1, \ldots , \mu_n)$, with $\mu_i$ the number of $i$'s in $T$.
\item $T$ has no violation entries.
\item The left-justified form of $\widehat T$ (the dual of $T$) is a SSYT or, equivalently, $\widehat{T}$ is highest weight (as a ptableau).
\end{enumerate} \label{word equiv}
\end{cor}

\begin{proof} The equivalence of (1) and (2) was proved in~\cite{ptab}, Lemma 7.6, and the equivalence of (1) and (3) is easy to show (note that the left-most violation entry in the right-justified form of $T$, if it exists, would show $T$ does not satisfy the word condition). The equivalence of (1) and (4) is proved by noting, from Corollary~\ref{blank dual cor} above, that $\widehat{T}^*$ has no blanks above content, and so by Theorem~\ref{highest weight ptab}, it is highest weight.

\end{proof}

For example, in the right-justified form:

\[ T=\begin{ytableau} \ &&&&&1&1&1\\&&&&1&&2&2\\&1&1&1&&2&&\\1&2&2&2&2&3&3&3\\&&&&&4&4&4\end{ytableau} \, , \]
we can see that $T$ satisfies the word condition by noting that every number $\ell$  in  $T$, for $2 \leq \ell$, has an $\ell -1$ above it.  $\widehat{T}$ is then partition-shaped (and hence highest weight as a ptableau):

\[ \widehat{T}= \begin{ytableau}1 & 1 & 1 & 2 & 3 & 3 & 3 & 4\\2&2&3&4&4&4&4&\\4&4&4&&&&&\\5&5&5&&&&&\end{ytableau}\, . \]

\section{Growth and Highest Weight Ptableaux}

The following result is an elementary observation, but the authors were unable to find it stated in the literature. It shows, among other things, that if a word $\omega=\omega_1 \omega_2 \cdots \omega_k$ is a node in a crystal of highest weight $\tau = \tau_1 \tau_2 \cdots \tau_k$, then the highest weight of the crystal of $\omega \otimes \omega_{k+1}$ is an extension of $\tau$, $\tau_1\tau_2\cdots\tau_k\tau_{k+1}'$. We express this result in the framework of ptableaux:

\begin{thm} Suppose $T$ is a perforated tableau. Let $T' = T\otimes \framebox{s}_{\,a}$ denote the perforated tableau produced by adjoining one box with value $s$ in row $a$ (see Definition~\ref{bi-word def}), such that the extension produces a valid ptableau. Then \label{epath high weight}
\begin{enumerate}
\item Given any sequence of crystal operators $e^{*} = e_{i_s} \cdots e_{i_1}$ such that $e^{*}T$ is highest weight, there is a sequence $e^{\star\star}$ determined only by $e^{*}$, the values $\phi_{i_{\ell}}(e_{i_{\ell-1}} \cdots e_{i_1} T)$, and the row index $a$, such that $e^{**}T'$ is highest weight. In particular, $e^{**}$ is independent of the value $s$, except that we require $T\otimes \framebox{s}_{\,a}$ to be a valid ptableau.
\item The highest weight ptableau $e^{**}T'$ extends the highest weight ptableau $e^{*}T$. That is, $e^{**}T'$ has the form of the highest weight ptableau $e^{*}T$ to which $\framebox{s}$ has been appended, in some row, to produce a partition-shaped ptableau.
\end{enumerate} \label{weight order}
\end{thm}

\begin{proof} Suppose the sequence $e^{*} = e_{i_s} \cdots e_{i_1}$ takes $T$ to its highest weight $e_{i_s} \cdots e_{i_1}T$.
We prove the existence of a modified sequence $e^{**}=e^{*i_s +1}e^{*i_s} \cdots e^{*i_1}$ such that $e^{**} T\otimes \framebox{s}_{\,a}$ is highest weight. In particular, we show:
\begin{description}
\item[Property 1.] Each $e^{*\ell}$ is a sequence of crystal operators determined only by $e^{*}$, the values $\phi_{i_{\ell}}(e_{i_{\ell-1}} \cdots e_{i_1} T)$, and the row index $a$, the row in which the box $\framebox{s}$ of $T'$ was first appended.
\item[Property 2.] For each $\ell \in \{ 1, \ldots, s \}$, the ptableau $e^{*\ell} e^{*(\ell -1)} \cdots e^{*1}T'$ is the ptableau $e_{i_{\ell}} \cdots e_{i_1} T$ to which the number $s$ has been appended, but possibly in a row above $a$.
\item[Property 3.] The ptableau $e^{**}T'$ is highest weight.
\end{description}

We define $e^{*i_\ell}$, inductively, as follows:
\begin{enumerate}
\item[a.] If $\phi_{*i_\ell}(e_{i_{\ell -1}}\cdots e_{i_1}T)=0$ and $\framebox{s}_{\,a}$ appears in row $i_\ell+1$ of $e^{*\ell} e^{*(\ell -1)} \cdots e^{*1}T'$, set $e^{*i_\ell}= e_{i_\ell} e_{i_\ell}$.
\item[b.] Otherwise, set $e^{*i_{\ell}}=  e_\ell$.
\end{enumerate}

Clearly Property 1 has been satisfied. We prove Property 2 now.

Consider Case (a): Suppose $\phi_{i_1}(T) =0$ and $\framebox{s}_{\,a}$ appears in row $i_1 +1$.
Recall that $\phi_{i_1} (T)$ equals the number of blanks in row ${i_1}+1$ of $^*T[{i_1},{i_1}+1]$, which we assume equals zero in this case. Thus:
\[^*T[{i_1},{i_1}+1] = \begin{ytableau} \none & \none & \none & c' & \cdots & \cdots & \none & \none & c''' & \cdots \\ \cdots & \cdots & c & a & \cdots & \cdots & \cdots & c''& b& \cdots  \end{ytableau} \]
but then
\[^*T'[{i_1},{i_1}+1] = \begin{ytableau} \none & \none & \none & c' & \cdots & \cdots & \none & \none & c''' & \cdots \\ \cdots & \cdots & c & a & \cdots & \cdots & \cdots & c''& b& \cdots & s \end{ytableau}. \]
We compute $e_{i_1}(T')$ by moving up the rightmost uncovered content, so that
\[^*(e_{i_1} T')[{i_1},{i_1}+1] = \begin{ytableau} \none & \none & \none & c' & \cdots & \cdots & \none & \none & c''' & \cdots&s \\ \cdots & \cdots & c & a & \cdots & \cdots & \cdots & c''& b& \cdots  \end{ytableau}. \]
Then, applying $e_{i_1}$ one additional time moves the \emph{same} entry that was moved when $e_{i_1}$ was applied to $^*T[{i_1},{i_1}+1]$, so we obtain:
\[^*(e_{i_1} e_{i_1} T')[i,i+1] = \begin{ytableau} \none & \none & \none & c' & \cdots & \cdots & \none & c'' & c''' & \cdots&s \\ \cdots & \cdots & c & a & \cdots & \cdots & \cdots & & b& \cdots  \end{ytableau}. \]
Setting $e^{*1} = e_{i_1} e_{i_1}$ satisfies Property 2.

For Case (b) above, suppose $\phi_{i_1}(T) > 0$, which means the left-justified $^*T[{i_1},{i_1}+1]$ has a ``shelf" in row ${i_1}$ jutting to the right, with blanks below:
\[^*T[{i_1},{i_1}+1] =  \begin{ytableau} \none & \none & \none & c' & \cdots & \cdots & \none & \none & c''' & \cdots & q& \cdots&r\\  \cdots & c & a & \cdots & \cdots & \cdots & c''& b& \cdots& p \end{ytableau}. \]
Suppose for now that $\framebox{s}$ appears in row ${i_1}+1$. Since $T'$ is a valid ptableau, we know that all entries in row ${i_1}$ must be strictly smaller than $s$, so that
\[^*T'[{i_1},{i_1}+1] =  \begin{ytableau} \none & \none & \none & c' & \cdots & \cdots & \none & \none & c''' & \cdots & q& \cdots&r\\  \cdots & c & a & \cdots & \cdots & \cdots & c''& b& \cdots& p&s \end{ytableau}. \]
Thus, applying $e_{i_1}$ to $T'$ moves the same right-most uncovered entry in $T'$ as applying $e_{i_1}$ to $T$ would move, so that setting $e^{*1}=e_{i_1}$ satisfies Property 2.

In the case that $\framebox{s}$ is appended in row ${i_1}$, then
\[^*T'[{i_1},{i_1}+1] =  \begin{ytableau} \none & \none & \none & c' & \cdots & \cdots & \none & \none & c''' & \cdots & q& \cdots&r&s\\  \cdots & c & a & \cdots & \cdots & \cdots & c''& b& \cdots& p \end{ytableau}, \]
and we see that setting $e^{*1}=e_{i_1}$ satisfies Property 2.\ as well.

To show Property 2 is implied by the subsequent $e^{*i_{\ell}}$ we apply the same reasoning as above, replacing $T$ with $e_{i_1}T$, then replacing $T$ with $e_{i_2}e_{i_1}T$, etc.

Thus, applying $e^{*i_s}\cdots e^{*i_1}$ to $T'$ produces a highest weight ptableau $e^{*}T$ to which a box $\framebox{s}$ has been appended in some row. In the case that $e^{*i_s}\cdots e^{*i_1}T'$ is not itself highest weight, then there must be some row in which $\epsilon_{i^*}(e^{*i_s}\cdots e^{*i_1}T') >0$, even though all $\epsilon_{i}(e_{i_s}\cdots e_{i_1}T) =0$. If so, we set $e^{*i_{\ell +1}}=\cdots e_{i* -1} e_{i^*}$ where we apply these operators until $e^{*i_{\ell +1}}(e^{*i_s}\cdots e^{*i_1}T')$ is highest weight. For example, if some portion of $(e^{*i_s}\cdots e^{*i_1}T')$, starting with row $i^* +1$ on the bottom up to some row $i^**$ at the top, had the shape:
\[\begin{ytableau} a & \cdots & \cdots & \cdots & \cdots & \cdots &a' & \cdots & \cdots & \cdots &a''\\
b &  \cdots & \cdots & \cdots & \cdots & \cdots &b'\\
\vdots & \cdots & \cdots & \cdots & \cdots & \cdots &\vdots \\
c& \cdots & \cdots & \cdots & \cdots & \cdots &c'&s\end{ytableau} \]
then applying $e_{i** +1} \cdots e_{i*}T'$ produces:
which now, being partition shaped, is highest weight. This last observation proves Property 3, completing the proof.
\end{proof}

For example, let $\Hw(T)$ denote the highest weight element in an irreducible crystal graph containing $T$. Then if
\[ T = \begin{ytableau} \ & 1 & 1 & 2 &3 &4 \\ &2&2&3&& \\ 1 &3&3&4&5&5\end{ytableau}\, , \ \ \hbox{we have} \ \ \Hw(T) = \begin{ytableau} 1 & 1& 1&2 &3&4 \\ 2 &2&3&5&5& \\ 3&3&4&& & \end{ytableau} \, . \]
But if we form the extension $T \otimes \framebox{$5$}_{\, 2}$, we see
\[ T\otimes \framebox{$5$}_{\, 2} = \begin{ytableau} \ & 1 & 1 & 2 &3 &4 &\\ &2&2&3&&&5 \\ 1 &3&3&4&5&5&\end{ytableau}\, , \ \ \hbox{and hence} \ \ \Hw(T\otimes \framebox{$5$}_{\, 2}) = \begin{ytableau} 1 & 1& 1&2 &3&4&5 \\ 2 &2&3&5&5& &\\ 3&3&4&& && \end{ytableau} \,. \]
That is, $\Hw(T\otimes \framebox{$5$}_{\, 2}) = \Hw(T) \otimes \framebox{$5$}_{\, 1}$, as predicted in the theorem.

\medskip
By restricting to the crystal structure of the bottom row $\omega$ of a biword $\tau \choose \omega$, we immediately obtain the following:

\begin{cor} Let $\omega=\omega_1\omega_2\cdots\omega_k  \in [n]^{\otimes k}$. Let $\omega'$ be some initial factor of $\omega$ (that is, $\omega' =\omega_1\omega_2\cdots \omega_j$ for some $j < k$). Then the highest weight of the irreducible crystal containing $\omega'$ is an initial factor of the highest weight of the irreducible crystal containing $\omega$.
\end{cor}

For example, using $T$ as in the example above, we see
\[ \Bw(T) = {\tau \choose \omega} = { 11122233334455 \choose 31122133213133  },\  \Hw{\tau \choose \omega} = {11122233334455 \choose 11122133213122} \, .\]
But then
\[ \Hw\left( {\tau  \choose \omega} \otimes { 5 \choose  2} \right) ={111222333344555 \choose 111221332131221} =  {11122233334455 \choose 11122133213122} \otimes {5 \choose 1} \]
and, in particular, $\Hw(\omega \otimes 2)$ = $\Hw(\omega) \otimes 1$.



\section{RSK and PTableaux}

Below we define a ptableau RSK insertion procedure that takes a ptableau $T$, and produces a pair of ptableaux $(PT,T_{\mmax})$ such that:

\begin{enumerate}
\item $T_{\mmax}$ is the highest weight element of the irreducible crystal containing $T$.
\item $PT$ is a ptableau plactically equivalent to $T$, but lying in the \emph{distinguished crystal} (recall Definition~\ref{distinguished crystal def}).
\end{enumerate}
Ê
\noindent As an immediate consequence of these facts, we will conclude:
\begin{itemize}
\item Given two ptableaux $T,T'$, if $RSK(T) = (PT,T_{\mmax})$ and $RSK(T') = (PT',T_{\mmax}')$, then $T$ and $T'$ lie in the same crystal only if $T_{\mmax} = T_{\mmax}'$.
\item Given two ptableaux $T,T'$, if $RSK(T) = (PT,T_{\mmax})$ and $RSK(T') = (PT',T_{\mmax}')$, then $T$ and $T'$ are plactically equivalent only if $PT = PT'$.
\end{itemize}

This approach to plactic equivalence is more natural than the classical RSK correspondence. To decide if two nodes of a crystal graph lie in the same connected component, instead of computing some \emph{new} object as a proxy for the highest weight nodes, we use an insertion scheme to compute highest weight elements directly. To tell if two nodes are plactically equivalent, we compute plactically equivalent representatives in the distinguished crystal of ptableaux. We describe our algorithm below, and prove these assertions. Later we establish the equivalence (duality) of our approach and the classical RSK algorithm.

\subsection{RSK Insertion for Ptableaux}
\begin{thm} The algorithm {\emph {RSK}}  (defined below) maps a ptableau $T$ to a pair of ptableaux $RSK(T)= (PT, T_{\mmax})$ such that
\begin{enumerate}
\item $PT$ is an element in the distinguished crystal.
\item $PT$ is plactically equivalent to $T$.
\item $T_{\mmax}$ is the highest weight element in the irreducible crystal graph containing $T$.
\end{enumerate} \label{RSK algorithm}
\end{thm}

\begin{proof} We first present the algorithm:

\bigskip

\noindent {\bf\texttt{Algorithm} {\emph {RSK}} }

Input: A ptableau $T \in \PTab_{(m,n)}^k$, or, equivalently, its biword $\bw(T) = {\tau \choose \omega}\in \Bw([m]\times [n])^{\otimes k}$.

Output: A pair of ptableaux, $PT$ and $T_{\mmax}$, which we denote by {\emph {RSK}}$(T) = (PT,T_{\mmax})$.
\medskip

Recall the notation from Definition~\ref{bi-word def} for ${\tau \choose \omega}^{(\ell)}$ and $T^{(\ell)}$.  We let $PT^{(0)}$ and $T_{\mmax}^{(0)}$ denote empty ptableaux. We construct $PT$ and $T_{\mmax}$ sequentially, by defining $PT^{(\ell)}$ and $T_{\mmax}^{(\ell)}$ for $1 \leq \ell \leq k$.

We simplify notation as follows. Let $\alpha_{\ell}$ denote the number of columns in $PT^{(\ell)}$. Let ${\cal C}_{\alpha_{\ell}-j+1}$ denote column $\alpha_{\ell}-j+1$ of $PT^{(\ell)}$. Fixing some $\ell$, for $0 \leq \ell <k$, let $\framebox{$a_1$}_{\, b_1} =\framebox{$\tau_{\ell+1}$}_{\,\omega_{\ell +1}}$. That is, we set $a_1 = \tau_{\ell+1}$ and $b_1 = \omega_{\ell+1}$. We construct $PT^{(\ell+1)}$ from $PT^{(\ell)} \otimes \framebox{$a_1$}_{\, b_1} $ by defining successive ``column insertions'' of some cell $\framebox{$a_j$}_{\, b_j} $ into the column ${\cal C}_{\alpha_{\ell}-j+1}$ , followed by right-justification and a step called \emph{resolution}, resulting in pushing a new cell to the \emph{left} of $C_{\alpha_{\ell}-j+1}$ which will be denoted $\framebox{$a_{j+1}$}_{\, b_{j+1}}$, unless the right-justification produces no such cell, and the insertion terminates.


\ytableausetup{smalltableaux}

It is simpler to prove assertion (1) (that the ptableau $PT$ we construct lies in the distinguished crystal) along with the presentation of the algorithm itself.
We therefore assume inductively that $PT^{(\ell)}$ satisfies the word condition, and demonstrate that the ptableau $PT^{(\ell+1)}$ that we produce also satisfies the word condition.
Thus, in some column ${\cal C}_{\alpha_{\ell}-j+1}$, we represent the insertion/resolution process by steps $(i), (ii)$ and $(iii)$ below:
\begin{figure}[H]
\begin{center}
{\includegraphics[scale=.27]{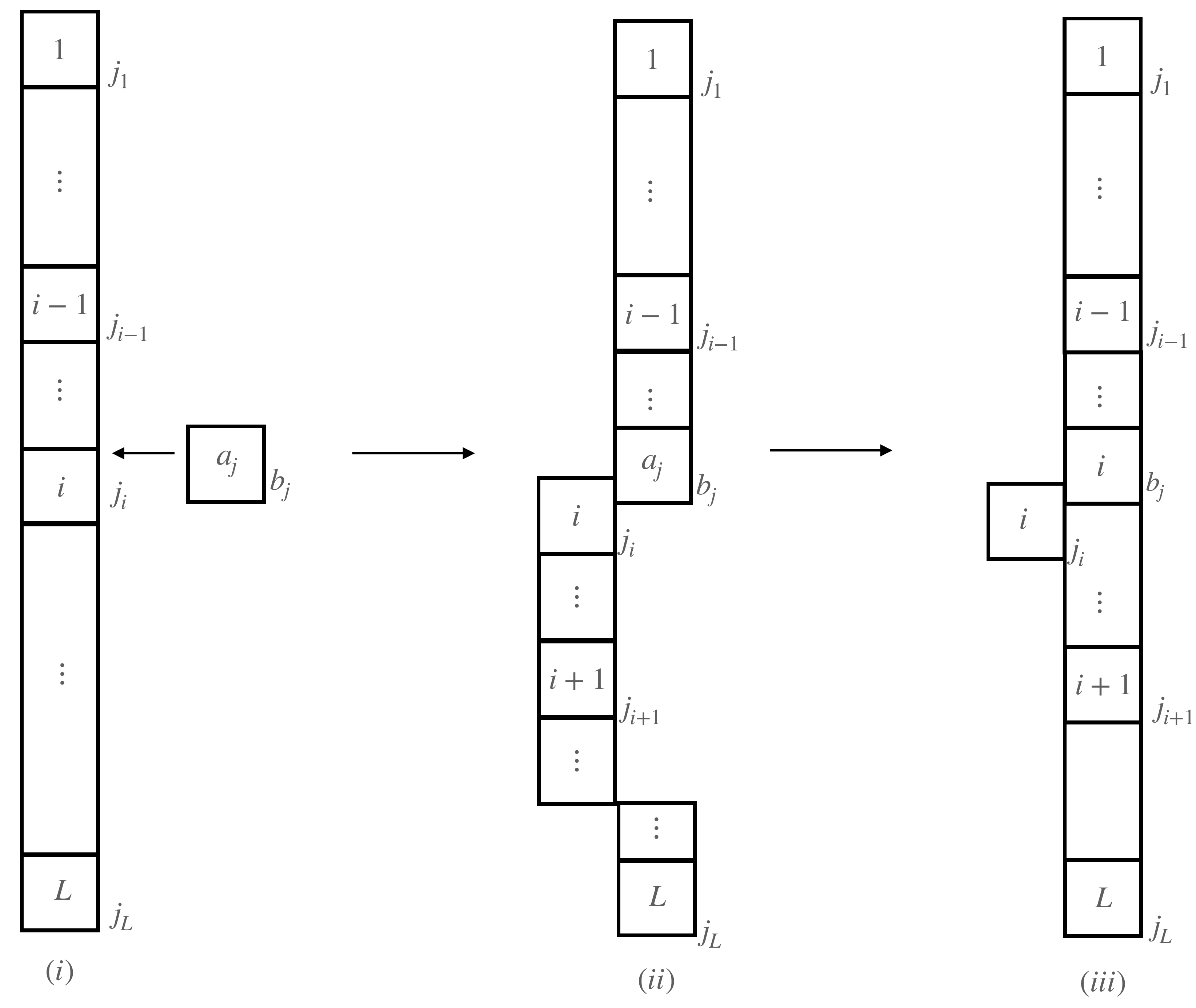}}
\end{center}
\end{figure}

The change from $(i)$ to $(ii)$ is due to the right-justification of the column with inserted content ${\cal C}_{\alpha_{\ell}-j+1}\otimes \framebox{$a_j$}_{\, b_j}$. The row index $j_i$ of the cell $\framebox{$i$}_{j_i}$ is determined as the minimum row index $j_{i'}$ such that $j_{i'} \geq b_j$. All cells $\framebox{$t$}_{j_t}$ at or below $\framebox{$i$}_{j_i}$ with content $t \leq a_j$ will also be pushed left as part of the right-justification.

The transition from step $(ii)$ to $(iii)$ is from resolution, whereby the content of the inserted cell $\framebox{$a_j$}_{\, b_j}$ is sequentially reduced until it equals $i$. We call \emph{each} decrement of the value $s$ of a cell to a value $s-1$ a \emph{swap} (we present an interpretation of each swap in the context of classical RSK insertions below). Once all swaps of the resolution process have been completed, all content previously below the cell $\framebox{$i$}_{j_i}$ will move under the resolved inserted cell $\framebox{$i$}_{\, k_j}$, except for the cell $\framebox{$i$}_{\, j_i}$, which will (inductively) form the inserted cell of column $\alpha_\ell-j$ to its left, so we set $\framebox{$a_{j+1}$}_{\, b_{j+1}}=\framebox{$i$}_{\, j_i}$ and continue left with the insertion into
${\cal C}_{\alpha_{\ell}-j}$.

Note that if $a_j > \ell$, then the inserted cell will push no content left, and so after resolution (when its content $s_j$ is decreased to $\ell +1$, if necessary), the insertion algorithm stops. We call this the \emph{terminal case}. In all cases, however, we see that by stage $(iii)$ of the resolution process, column $\alpha_\ell-j+1$ will have no violation entries for $j=1, \ldots $, and hence, upon reaching the terminal case, $PT^{(\ell+1)}$ will satisfy the word condition and hence, lie in the distinguished crystal, thus proving assertion (1).

\medskip

We complete the description of the algorithm. To compute the updated $T_{\mmax}^{(\ell+1)}$, we let $\eta_{\ell+1}$ denote the \emph{content} of the last inserted cell in the terminal case. We then set
\[ T_{\mmax}^{(\ell+1)} = T_{\mmax}^{(\ell)} \otimes \framebox{$\tau_{\ell+1}$}_{\,\eta_{\ell +1}}, \]
recalling that $\tau_{\ell+1}$ is the value of the content inserted in $T^{(\ell)}$ at this stage.
That is, the \emph{value} of the terminal insertion in the construction of $PT^{(\ell +1)}$ determines the \emph{row} in which we add the content $\framebox{$\tau_{\ell+1}$}$ to update $T_{\mmax}^{(\ell)}$ to obtain $T_{\mmax}^{(\ell+1)}$.

\medskip

Proof of (2): Recall that if two ptableaux $S$ and $T$ are plactically equivalent, we write $S\equiv T$.  We want to prove that $PT^{(\ell+1)} \equiv T^{(\ell+1)}$, under the hypothesis $PT^{(\ell)} \equiv T^{(\ell)}$. We note the base case (for $\ell = 1$) is trivially satisfied:
$$T^{(1)} =\Pf{\binom{\tau_1}{\omega_1}} = \framebox{$\tau_{1}$}_{\, \omega_1}
\hbox{ and }
 PT^{(1)} =\Pf{\binom{1}{\omega_1}}=\framebox{$1$}_{\, \omega_1},$$ which are clearly plactically equivalent ptableaux.

So, suppose $PT^{(\ell)} \equiv T^{(\ell)}$. To prove $PT^{(\ell+1)} \equiv T^{(\ell+1)}$, we first claim that
$PT^{(\ell)} \otimes \framebox{$\tau_{\ell+1}$} \equiv PT^{(\ell+1)}$.

To prove this claim, we prove a (slightly) stronger claim, namely, that the ptableau formed from $PT^{(\ell)} \otimes \framebox{$\tau_{\ell+1}$}$ after \emph{each} swap of the resolution produces a ptableau that is plactically equivalent to the one proceeding it in the process. By the chain of plactic equivalences from $PT^{(\ell)} \otimes \framebox{$\tau_{\ell+1}$}$ to each of the ptableaux obtained by swapping, to the final ptableau $PT^{(\ell+1)}$, we prove the claim.

We argue that the swapping step of the resolution process commutes with the application of any crystal operator $f_i$. From this, the plactic equivalence claim follows.

Inductively, therefore, we consider applying a crystal operator $f_i$ to a ptableau $PT^{(\bigtriangleup)}$ that satisfies the word condition \emph{except} for the possibility of some entries (violation entries) in the right-justified $PT^{(\bigtriangleup)*}$, and we let $\framebox{$a$}_{\,t}$ denote the \emph{right-most} such violation entry. Thus, the next stage of resolution would be to swap the content of $\framebox{$a$}_{\,t}$ to $\framebox{$a-1$}_{\,t}$.

Suppose we consider applying some non-NULL crystal operator $f_i$ to $PT^{(\bigtriangleup)}$, where the box $\framebox{$c$}_{\,i}$ is moved by the application of $f_{i}$, so that this box becomes $\framebox{$c$}_{\,i+1}$ in the ptableau $f_i PT^{(\bigtriangleup)}$. Though there are, in principle, several cases to consider, they can be reduced to:

\medskip
\noindent Case 1: The box $\framebox{$a$}_{\,t}$ is in a \emph{different column} than $\framebox{$c$}_{\,i}$ in the right-justified form of $PT^{(\bigtriangleup)}$.

\medskip
\noindent Case 2: The boxes $\framebox{$a$}_{\,t}$ and $\framebox{$c$}_{\,i}$ are in the \emph{same column} in $PT^{(\bigtriangleup)}$.

\bigskip

\noindent Proof for Case 1: There are several potential cases to consider, but they all amount to the same observation that in the right-justified $PT^{(\bigtriangleup)*}$, applying the crystal operator $f_i$ cannot change the entry in $PT^{(\bigtriangleup)*}$ that we must swap during resolution, nor can applying the swap alter the box that is moved by $f_i$. For example, suppose the to-be-swapped box $\framebox{$a$}_{\,t}$ is in some column that is strictly to the \emph{left} of the column containing the entry $\framebox{$c$}_{\,i}$ that will be moved by application of $f_i$:

\[ \begin{ytableau}
\none[\hbox{Row $i$:\ \ }]&\none & \none& \scriptstyle{\dots} & a&\scriptstyle{\dots}&c&\scriptstyle{\dots} \\
\none[\hbox{Row $i+1$:\ \ \  }] &\none &\none& \scriptstyle{\dots}& \star & \scriptstyle{\dots}& & \scriptstyle{\dots} \end{ytableau}. \]
Note that, in the diagram above, we are also assuming that $\framebox{$a$}$ is \emph{also} in row $i$, although that will not be necessary for our analysis. Here, the ``$\dots$" denote the right-justified content in rows $i$ and $i+1$ of $PT^{(\bigtriangleup)*}$. There is a blank under $\framebox{$c$}$ that will be moved by the application of $f_i$, and the $\star$ must denote some content strictly larger than $a$ (else the $f_i$ moved box would not be $\framebox{$c$}$\,). From this, it then follows that $a < \star \leq c$. We claim is that the application of $f_i$ and the swapping operation on $PT^{(\bigtriangleup)*}$ commute:

\bigskip
\ytableausetup{nosmalltableaux}

\[ \begin{array}{ccc}
\begin{ytableau}
 \scriptstyle{\dots} & a&\scriptstyle{\dots}&c&\scriptstyle{\dots} \\ \scriptstyle{\dots}& \star & \scriptstyle{\dots}& & \scriptstyle{\dots} \end{ytableau} & \stackrel{f_i}{\longrightarrow} & \begin{ytableau}
\scriptstyle{\dots} & a&\scriptstyle{\dots}&&\scriptstyle{\dots} \\ \scriptstyle{\dots}& \star & \scriptstyle{\dots}&c & \scriptstyle{\dots} \end{ytableau} \\
\downarrow \hbox{Swap} & & \downarrow \hbox{Swap} \\
\begin{ytableau}
 \scriptstyle{\dots} & \scalebox{0.5}[0.7]{$a-1$}&\scriptstyle{\dots}&c&\scriptstyle{\dots} \\ \scriptstyle{\dots}& \star & \scriptstyle{\dots}& & \scriptstyle{\dots} \end{ytableau} & \stackrel{f_i}{\longrightarrow} & \begin{ytableau}
 \scriptstyle{\dots} & \scalebox{0.5}[0.7]{$a-1$}&\scriptstyle{\dots}&&\scriptstyle{\dots} \\ \scriptstyle{\dots}& \star & \scriptstyle{\dots}&c & \scriptstyle{\dots} \end{ytableau} \end{array} \]
 
 \bigskip

 The commutativity follows by noting that the relation $a < \star \leq c$ implies that applying $f_i$ will not remove the violation entry at $\framebox{$a$}$ (since moving content with a \emph{higher} value than $a$ cannot affect a word condition involving the number of $a-1$'s over the $a$'s) and, similarly, swapping the value of $a$ for $a-1$ will not change which box is moved by $f_i$ (as swapping the value does not change the location of blanks, so that the $f_i$-moved content remains $c$). As said above, there are several other such cases, but they all amount to the same observation, in that the $f_i$-moved box remains in the same column after applying $f_i$, but it cannot affect the need for swapping at the violation entry.

 \medskip
 \noindent Proof for Case 2: As before, there are several cases, but only one  is not obvious: when the $f_i$-moved box is \emph{also} the violation entry. In all other cases, the two operations act independently of each other. Suppose that the $f_i$-moved content and the swapped value are both $a$. Since, before insertion of a cell, all content satisfied the word condition, there can only be a \emph{single} violation entry in any row, and so to the \emph{left} of that violation entry, the values are \emph{less than} $a$. If these values are less than $a-1$, then the moved cell $\framebox{$a$}_{\, i}$ remains the swapped cell after (or before) the move, similar to Case (1). If the values to the left of $\framebox{$a$}_{\, i}$ are $a-1$, then, in rows $i$ and $i+1$ we might have:
\ytableausetup{nosmalltableaux}
 \[ \begin{ytableau}\none[\hbox{Row $i$:\ \ }]&\none & \none &\scriptstyle{\dots} &\scalebox{0.5}[0.7]{$a-1$}&
 \scriptstyle{\dots} & \scalebox{0.5}[0.7]{$a-1$}&a&\scriptstyle{\dots}& \scriptstyle{\dots}&\scriptstyle{\dots} \\ \none[\hbox{Row $i+1$:\ \ }]&\none & \none&\scriptstyle{\dots}&\star&\scriptstyle{\dots}& \star & &\scriptstyle{\dots}&\scriptstyle{\dots}& \scriptstyle{\dots} \end{ytableau}\, .\]
 The notation
 \raisebox{0.05in}{$\begin{ytableau}\scalebox{0.5}[0.7]{$a-1$}&
 \scriptstyle{\dots} & \scalebox{0.5}[0.7]{$a-1$}\end{ytableau}$}
 in row $i$ above denotes the \emph{maximal} consecutive string of content with value $a-1$ in row $i$. In this case, we first conclude that the right-most $\star$ in row $i+1$ must have value $a$, or else we could move the value further to the right:

  \[ \begin{ytableau}\none[\hbox{Row $i$:\ \ }]&\none & \none &\scriptstyle{\dots} &\scalebox{0.5}[0.7]{$a-1$}&
 \scriptstyle{\dots} & \scalebox{0.5}[0.7]{$a-1$}&a&\scriptstyle{\dots}& \scriptstyle{\dots}&\scriptstyle{\dots} \\ \none[\hbox{Row $i+1$:\ \ }]&\none & \none&\scriptstyle{\dots}&\star&\scriptstyle{\dots}& a & &\scriptstyle{\dots}&\scriptstyle{\dots}& \scriptstyle{\dots} \end{ytableau}\, .\]

 We then conclude that \emph{every} $a-1$ appearing in row $i$ must have an $a$ beneath it:
   \[ \begin{ytableau}\none[\hbox{Row $i$:\ \ }]&\none & \none &\scriptstyle{\dots}&\scriptstyle{\dots} &\scalebox{0.5}[0.7]{$a-1$}&
 \scriptstyle{\dots} & \scalebox{0.5}[0.7]{$a-1$}&a&\scriptstyle{\dots}\\ \none[\hbox{Row $i+1$:\ \ }]&\none & \none&\scriptstyle{\dots}&s&a&\scriptstyle{\dots}& a & &\scriptstyle{\dots} \end{ytableau} \, .\]
 This follows by noting that there must \emph{some} content under them, or else the left-most covered blank would not be under the $a$, and this entry would then \emph{not} be the $f_i$-moved box, as assumed. Further, all these non-blank entries must have content $a$, since any entry below $a-1$ is at least $a$ by column strictness, but no more than $a$ since these entries are to the left of $a$ in row $i+1$. Finally, we must have $s <a$.  Otherwise, if we did have $s=a$ and there were (by maximality) no $a-1$ over it and if the left-most $a-1$ in the row above were not above $s$, this would force the cell with content $s$ to be a violation entry \emph{prior} to the insertion of the $a$ in row $i$ to its right, contradicting our hypothesis.

 With these facts in hand, we then obtain the following commutative diagram (recalling that we right-justify after the swap):

 \[ \begin{array}{ccc}
\begin{ytableau}\scriptstyle{\dots} &\scriptstyle{\dots} &\scalebox{0.5}[0.7]{$a-1$}&
 \scriptstyle{\dots} & \scalebox{0.5}[0.7]{$a-1$}&a&\scriptstyle{\dots}\\ \scriptstyle{\dots}&s&a&\scriptstyle{\dots}& a & &\scriptstyle{\dots} \end{ytableau} & \stackrel{f_i}{\longrightarrow} & \begin{ytableau}\scriptstyle{\dots}&\scriptstyle{\dots}  &\scriptstyle{\dots} &
  \scalebox{0.5}[0.7]{$a-1$} & \scriptstyle{\dots} & \scalebox{0.5}[0.7]{$a-1$}&\scriptstyle{\dots}\\ \scriptstyle{\dots}  & s&a&\scriptstyle{\dots}& a & a &\scriptstyle{\dots} \end{ytableau} \\
\downarrow \hbox{Swap} & & \downarrow \hbox{Swap} \\
\begin{ytableau}\scriptstyle{\dots}  &\scriptstyle{\dots} &\scalebox{0.5}[0.7]{$a-1$}&
 \scriptstyle{\dots} & \scalebox{0.5}[0.7]{$a-1$}&\scalebox{0.5}[0.7]{$a-1$}&
 \scriptstyle{\dots}\\ \scriptstyle{\dots}&s&&a&
 \scriptstyle{\dots} & a&\scriptstyle{\dots} \end{ytableau} & \stackrel{f_i}{\longrightarrow} &\begin{ytableau}\scriptstyle{\dots}  &\scriptstyle{\dots} &\scriptstyle{\dots} &
  \scalebox{0.5}[0.7]{$a-1$} & \scriptstyle{\dots} & \scalebox{0.5}[0.7]{$a-1$}&\scriptstyle{\dots}\\ \scriptstyle{\dots}&s& \scalebox{0.5}[0.7]{$a-1$}&\scriptstyle{\dots}& a & a &\scriptstyle{\dots} \end{ytableau} \end{array}. \]

Following the swap on the left we see the string of $a$'s in row $i+1$ shift one box to the right, so that the $f_i$-moved content is now \emph{moved} left to the left-most $a-1$ in row $i$, which then moves down after applying the $f_i$ operator in the bottom row of the diagram. However, applying the $f_i$ operator \emph{first} in the top row of the diagram moves an $a$ into row $i+1$, so that the $a-1$'s in row $i$ shift one box to the right. Thus, the violation entry in row $i+1$ shifts from the left-most $a$ in row $i$ to the left-most $a$ in row $1+1$, so that after the swap in that row, the diagram is the same.

Thus, at each step the swapping dictated by resolution commutes with the application of the $f_{i}$ operators, and so $PT^{(\ell)} \otimes \framebox{$\tau_{\ell+1}$} \equiv PT^{(\ell+1)}$.

This claim will (finally) allow us to prove $T \equiv PT$.
 Let $T^{(k<)}$ denote the ptableau satisfying, for each $k$,
$T =  T^{(k)}\otimes T^{(k<)}$ (recall the notation for $T^{(k)}$ in Definition~\ref{bi-word def}).
We can extend plactically equivalent words by a common extension, preserving plactic equivalence, and hence the same holds for ptableaux. We note that $T^{(1)}= \framebox{$\tau_1$}_{\, \omega_1}$ is plactically equivalent to $PT^{(1)}=\framebox{$1$}_{\, \omega_1}$ by construction, so that we have:

\begin{align*}
T =  T^{(1)} \otimes T^{(1<)}& \equiv PT^{(1)} \otimes T^{(1<)}\\
&= \left(PT^{(1)} \otimes \framebox{$\tau_2$}_{\, \omega_2}\right) \otimes T^{(2<)}\\
& \equiv PT^{(2)}  \otimes T^{(2<)}\\
& \equiv PT^{(3)}  \otimes T^{(3<)}\\
&\equiv \cdots \equiv PT^{(k)}=
 PT,
\end{align*}
which proves statement (2) of the Theorem.

Proof of (3): Trivially, $T_{\mmax}^{(1)}$ is highest weight in the crystal containing $T^{(1)}$, and, inductively, the highest weight ptableau lying in the crystal containing $T^{(j)}$ is $T_{\mmax}^{(j)}$. Extending $T^{(j)}$ to $T^{(j+1)}$ we conclude, by Theorem~\ref{epath high weight}, the highest weight ptableau containing $T^{(j+1)}$ extends $T_{\mmax}^{(j)}$ by the last entry of $T^{(j+1)}$ appended to $T^{(j)}$. Thus, it is sufficient to show that this last entry appended to $T^{(j+1)}$ appears in the same location in the highest weight as it does in $T_{\mmax}^{(j+1)}$.

We have $T^{(j)}$ and $PT^{(j)}$ are plactically equivalent by Statement (2). Since $PT^{(j)}$ satisfies the word condition, if the resolution producing $PT^{(j+1)}$ from $PT^{(j)} \otimes \framebox{$\tau_{j+1}$}_{\, \omega_{j+1}}$ terminates with a last entry with value $\eta_{j+1}$ in $PT^{(j+1)}$, then that element will necessarily appear in row $\eta_{j+1}$ in the highest weight element in the crystal containing $PT^{(j+1)}$ (since the highest weight elements for ptableaux in the distinguished crystal have only $j$'s in row $j$). That is, in the transition from $PT^{(j)}$ to $PT^{(j+1)}$, their respective highest weight ptableaux differ by one entry in row $\eta_{j+1}$.  By (2) of Theorem~\ref{epath high weight}, the highest weight of $T^{(j+1)}$ extends the highest weight of $T^{(j)}$, and by plactic equivalence it, too, must extend in row $\eta_{j+1}$, and thus, the highest weight of the extension of the crystal containing $T_{\mmax}^{(j)}$ is $T_{\mmax}^{(j+1)}$, which proves the theorem.
\end{proof}

\ytableausetup{smalltableaux}
\[ \hbox{\emph{Example:} We compute $T \mapsto \RSK(T)=(PT, T_{\mmax})$ for the ptableau} \ T =\begin{ytableau} \ & & 1 & 3 & 4 \\  1 & 2 & 2 && \\ 3 & 3 & 4 & 4& \end{ytableau}\, . \]
\vspace{-.3in}
\[ \renewcommand{\arraystretch}{2.5}
\begin{array}{cccc}k & T^{(k)} & PT^{(k)} & T_{\mmax}^{(k)}  \\ \hline
1 & \begin{ytableau} \ \\ 1 \\ \\ \end{ytableau} &  \begin{ytableau} \ \\ 1 \\ \\ \end{ytableau} & \begin{ytableau} 1 \\ \\ \\ \end{ytableau} \\ \hline
2 & \begin{ytableau} \ & 1 \\ 1 & \\& \\  \end{ytableau} & \begin{ytableau} \ & 1 \\ 1 & \\ &\\  \end{ytableau} & \begin{ytableau} 1 & 1 \\& \\ &\\ \end{ytableau} \\ \hline
3& \begin{ytableau} \ & 1 \\ 1 & 2 \\& \\ \end{ytableau} & \begin{ytableau} \ & 1 \\ 1 & 2 \\& \\ \end{ytableau} & \begin{ytableau} 1 & 1 \\ 2 & \\ &\\ \end{ytableau} \\ \hline
4 & \begin{ytableau} \ & & 1 \\ 1 & 2 &2 \\ & & \\ \end{ytableau}&
\begin{ytableau} \ & & 1 \\ 1 & 1 &2 \\ & & \\ \end{ytableau} & \begin{ytableau} 1 & 1 & 2 \\ 2 & & \\ && \\ \end{ytableau} \\ \hline
5 & \begin{ytableau} \ & & 1 \\ 1 & 2 & 2 \\ & & 3 \\ \end{ytableau}
& \begin{ytableau} \ & & 1 \\ 1 & 1 & 2 \\ & & 3 \\ \end{ytableau} &
\begin{ytableau} 1 & 1 & 2 \\ 2 & & \\ 3 & & \\ \end{ytableau}\\ \hline
6 & \begin{ytableau} \ & & 1 \\ 1 & 2 & 2 \\ & 3 &3 \end{ytableau} &
 \begin{ytableau} \ & & 1 \\ 1 & 1 & 2 \\ & 2 &3 \end{ytableau} &
 \begin{ytableau} 1 & 1 & 2 \\ 2 & 3 & \\ 3 & & \end{ytableau} \\ \hline
7& \begin{ytableau} \ & & 1 &3 \\ 1 & 2 & 2 & \\ & 3 & 3 & \end{ytableau} &
 \begin{ytableau} \ & & 1 &1 \\ 1 & 1 &  &2 \\ & &2 & 3   \end{ytableau} &
 \begin{ytableau} 1 & 1 & 2 & 3 \\ 2 & 3 & & \\ 3 &&& \end{ytableau} \\ \hline
 8& \begin{ytableau} \ &&1 & 3 \\ 1 & 2 & 2 & \\ & 3 & 3 & 4 \end{ytableau} &\begin{ytableau} \ &&1 & 1 \\ 1 & 1 &  &2 \\ & 2 & 2 & 3 \end{ytableau} &
 \begin{ytableau} 1 & 1 & 2 & 3 \\ 2 & 3 & 4 & \\ 3 & & & \end{ytableau} \\ \hline
9& \begin{ytableau} \ & & 1 & 3 \\ 1 & 2 & 2 &\\3 & 3& 4 & 4 \end{ytableau} & \begin{ytableau} \ & & 1 & 1 \\ 1 & 1 &  &2\\ 2& 2& 2 & 3 \end{ytableau} &
\begin{ytableau} 1 & 1 & 2 & 3 \\ 2 & 3 & 4 & 4\\ 3&&& \end{ytableau} \\ \hline
10 & T=\begin{ytableau} \ & & 1 & 3 & 4 \\  1 & 2 & 2 & &\\ 3 & 3 & 4 & 4& \end{ytableau} &
PT=\begin{ytableau} \ & & 1 & 1 & 1 \\ 1& 1 &  &  & 2\\ &2 & 2 & 2 & 3 \end{ytableau} &
T_{\mmax}=\begin{ytableau} 1 & 1 & 2 & 3 & 4 \\ 2 & 3 & 4 & 4& \\ 3 &&&& \end{ytableau} \\ \hline
\end{array}\]

With the ptableau RSK algorithm given and its properties proved, it's worth noting that although we defined the insertion column by column, in actuality, we can obtain $PT$ much more simply. In particular, resolving the insertion $PT^{(k)}\otimes \framebox{$\tau_{k+1}$}_{\, \omega_{k+1}}$ can be completed by resolving the violation entries of the right-justified form directly, without needing to resort to column-by-column insertion.

For example, if
\[ PT^{(k)}=\begin{ytableau} \ &&&&&1&1&1\\&&&&1&&2&2\\&1&1&1&&2&&\\1&2&2&2&2&3&3&3\\&&&&&4&4&4\end{ytableau} \, , \]
then resolving $PT^{(k)}\otimes \framebox{$4$}_{\,3}$ amounts to changing the violation entries (circled below) from right to left:

\ytableausetup{smalltableaux}
\[ \begin{ytableau} \ &&&&&&1&1&1\\&&&&&1&&2&2\\&&1&1&1&&2&&\bf \textcircled{$\scriptstyle \bf 4$}&\none\\  1& \textcircled{$\scriptstyle \bf 2$}& 2& 2& 2& \textcircled{$\scriptstyle \bf 3$}& 3& 3&\\&&&&& 4& 4& 4&\end{ytableau} \Rightarrow
\begin{ytableau} \ &&&&&&1&1&1\\&&&&&1&&2&2\\&&1&1&1&&2&&\bf \textcircled{$\scriptstyle  3$}&\none\\  1& \textcircled{$\scriptstyle \bf 2$}& 2& 2& 2& \textcircled{$\scriptstyle \bf 3$}& 3& 3&\\&&&&&& 4& 4& 4\end{ytableau}\Rightarrow
\begin{ytableau} \ &&&&&&1&1&1\\&&&&&1&&2&2\\&&1&1&1&&2&&\bf \textcircled{$\scriptstyle  3$}&\none\\  1& \textcircled{$\scriptstyle \bf 2$}& 2&\bf 2& 2& \textcircled{$\scriptstyle  2$}& 3& 3&\\&&&&&& 4& 4& 4\end{ytableau}\Rightarrow
\begin{ytableau} \ &&&&&&1&1&1\\&&&&&1&&2&2\\&&1&1&1&&2&&\bf \textcircled{$\scriptstyle  3$}&\none\\  1& \textcircled{$\scriptstyle 1$}& 2& 2& 2& \textcircled{$\scriptstyle  2$}& 3& 3&\\&&&&&& 4& 4& 4\end{ytableau}. \]

Notice that initially changing the 4 to 3 in the rightmost column resolves a word condition violation of 2's over 3's in rows 2 and 3, but then creates a new violation of 2's and 3's in rows 2 and 3 over rows 3 and 4. This is resolved in the second step, but then creates a word violation of 1's over 2's in rows 1 through 3 over rows 2 through 4, which is finally rectified with the final resolution step.
\begin{lem}
Suppose $PT$ satisfies the word condition, and let $P= \widehat{PT}$ (see Definition~\ref{dual def}). Let $PT'$ denote the resolution of $PT \otimes \framebox{$s$}_{\, k}$. Then $\widehat{PT'}$ equals the SSYT obtained by RSK column insertion of $k$ into $P$.

In particular, each swap of the resolution of $PT \otimes \framebox{$s$}_{\, k}$ corresponds, dually, with applying a crystal operator $e_i$ to $\widehat{PT \otimes \framebox{$s$}_{\, k}}$, taking it to a highest weight ptableau dual to $PT^{*}$.

\label{insert RSK}
\end{lem}
\begin{proof} We show that resolving each column (right to left) in $PT \otimes \framebox{$s$}_{\, k}$ is dual to the corresponding column (left to right) produced by RSK column insertion of $k$ into $P$, where $P= \widehat{PT}$. Let $\alpha$ denote the number of columns in $PT \otimes \framebox{$s$}_{\, k}$. Then, as in Theorem~\ref{RSK algorithm}, let ${\cal C}_{\alpha-j+1}$ denote column $\alpha-j+1$ of $PT \otimes \framebox{$s$}_{\, k}$. Let ${\cal C}_{\alpha-j+1}^*$ denote the column obtained after resolution of the entry in the column, where we assume, working right to left, that columns $\alpha$ through $\alpha-j+2$ have been resolved. It is sufficient to note that the columns $(i), (ii)$, and $(iii)$ of the insertion and resolution of ${\cal C}_{\alpha-j+1}\otimes \framebox{$s_j$}_{\, k_j}$ in column $\alpha-j+1$ of the resolution process (see the diagram in the proof of (1), Theorem~\ref{RSK algorithm}) are, in fact, dual to $(i), (ii)$ and $(iii)$ of column $j$ in the dual $\widehat{PT \otimes \framebox{$s$}_{\,k}}$:

\begin{figure}[H]
\begin{center}
{\includegraphics[scale=.27]{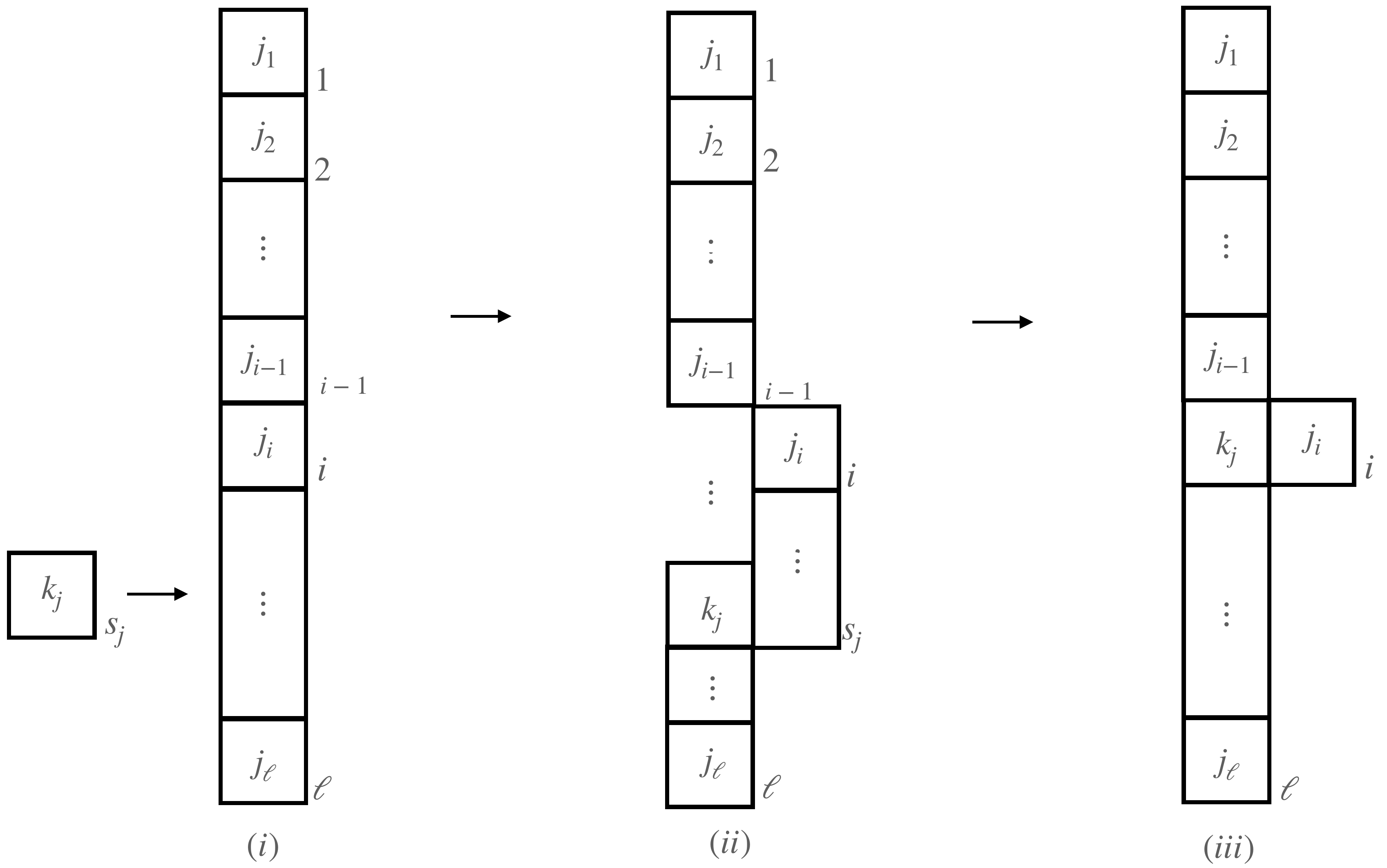}}
\end{center}
\end{figure}

That is, if there is a cell $\framebox{$i$}_{j_i}$ in column $\alpha-j+1$ of $PT$, then, reversing the role of row and content, there is a cell $\framebox{$j_k$}_{i}$ in column $j$ of $\widehat{PT}$, according to Theorem~\ref{left-right}. In particular, this column, before the (left) insertion of $\framebox{$k_j$}_{s_{j}}$, has content in consecutive rows $1$ through $\ell$.

Thus, upon \emph{left} justification of the content in $(i)$, we obtain $(ii)$, in which the cells in rows $s_j$ and higher with content $j_{i'}$ such that $k_j \leq j_{i'}$ are pushed right upon left justification. The resolution in ${\cal C}_{\alpha-j+1}\otimes \framebox{$s_j$}_{\, k_j}$ changing $\framebox{$s_j$}_{\, k_j}$ to $\framebox{$i$}_{\, k_j}$ has the dual effect, by virtue of Corollary~\ref{blank dual cor}, of moving the cell $\framebox{$k_j$}_{s_{j}}$ up into row $i$, after which all cells in rows $i+1$ and lower return left to column $j$ again, so that the cell $\framebox{$j_i$}_{i}$ is pushed right, forming the inserted content into column $j+1$ to its right. Since the content at a violation entry \emph{decreases} with resolution, and the content at the next violation entry to the right decreases again as we move right to left in $PT \otimes \framebox{$s$}_{\, k}$, we see moving the cell $\framebox{$s_j$}_{\, k_j}$ upwards in the dual $\widehat{PT \otimes \framebox{$s$}_{\,k}}$ amounts to applying a crystal operator $e_i$ in which there is only \emph{one} cell in the corresponding row with a blank over it.

In particular, we see that the transition $(i) \rightarrow (iii)$ in the above diagram is \emph{precisely} RSK column insertion of the content $k_j$ in column $j$, in which the cell with content $j_i$ is chosen as the highest such that $k_j \leq j_{i}$, and it becomes the content ``bumped'' into the column to its right.

Note, finally, that if the insertion ${\cal C}_{\alpha-j+1}\otimes \framebox{$s_j$}_{\, k_j}$ is the terminal case, then dually, in column $j$ we are inserting $\framebox{$k_j$}_{s_j}$ in a row \emph{lower} than $\ell$, so that, after resolution, the content will appear in row $\ell +1$, so that our insertion and resolution in ${\cal C}_{\alpha-j+1}\otimes \framebox{$s_j$}_{\, k_j}$ terminates at precisely the step that, dually, RSK column insertion would terminate with adding $\framebox{$k_j$}_{\, s_j}$ to the bottom of column $j$, and thus the lemma is proved.
\end{proof}

For example, if $PT^{(k)}$ were the ptableau below (so that $k=21$):
\ytableausetup{smalltableaux}
\[ PT^{(k)}=\begin{ytableau} \ &&&&&1&1&1\\&&&&1&&2&2\\&1&1&1&&2&&\\1&2&2&2&2&3&3&3\\&&&&&4&4&4\end{ytableau} \, ,
\ \ \hbox{its dual is} \ \ \widehat{PT^{(k)}}=
 \begin{ytableau} 1 & 1& 1&2&3&3&3&4\\2&2&3&4&4&4&4&\\4&4&4&&&&&\\5&5&5&&&&&\end{ytableau}\, ,\]
a highest weight ptableau (SSYT). But then
\[ PT^{(k)} \otimes \framebox{$4$}_{\, 3}= \begin{ytableau} \ &&&&&&1&1&1\\&&&&&1&&2&2\\&&1&1&1&&2&&\bf \textcircled{$\scriptstyle \bf 4$}&\none\\ \bf 1& \textcircled{$\scriptstyle \bf 2$}&\bf 2&\bf 2&\bf 2& \textcircled{$\scriptstyle \bf 3$}&\bf 3&\bf 3&\\&&&&&\bf 4&\bf 4&\bf 4&\end{ytableau} \ \ \hbox{and its dual is} \ \ \widehat{PT^{(k)} \otimes \framebox{$4$}_{\, 3}}=\begin{ytableau} 1 & 1& 1&2&3&3&3&&4\\2&2&3&&4&4&4&\textcircled{$\scriptstyle \bf 4$}&\\&4&4&\textcircled{$\scriptstyle \bf 4$}&&&&&\\\textcircled{$\scriptstyle \bf 3$}&5&5&5&&&&&\end{ytableau} \, , \]
where we have circled the content under blanks corresponding to the violation entries in $PT^{(k)} \otimes \framebox{$4$}_{\, 3}$, in reverse order to the blanks over content in its corresponding dual. Below we show the sequence $(i), (ii)$, and $(iii)$ for insertion into $PT$ from the proof in the first column (right justification and resolution), followed by ``left insertion'' into the dual $\widehat{PT}$ in the second column (left justification, with applying appropriate crystal operators $e_i$ (in the example, it is $e_3$) to return the ptableau to highest weight), and finally, classic RSK column insertion into $P= \widehat{PT}$, demonstrating that the output of left insertion into $\widehat{PT}$ equals (and explains) why RSK column insertion produces a highest with ptableau dual to $PT$.

\ytableausetup{smalltableaux}
\[ \begin{array}{rccc} & PT & \widehat{PT} & P\\
&&& \\

 & \begin{ytableau} \ &&&&&1&1&1\\&&&&1&&2&2\\&1&1&1&&2&& &\none[\leftarrow] & \none[\framebox{4}_{\,3}] \\1&2&2&2&2&3&3&3\\
&&&&&4&4&4\end{ytableau} &
\begin{ytableau}\none & \none&1 & 1& 1&2&3&3&3&4\\\none & \none &2&2&3&4&4&4&4&\\ \none & \none & 4&4&4&&&&&\\\none[\framebox{$3$}_{\, 4}] & \none[\rightarrow] & 5&5&5&&&&&\end{ytableau}
&
\begin{ytableau} \none & \none &1 & 1& 1&2&3&3&3&4\\\none & \none &2&2&3&4&4&4&4&\\ \none[3] & \none[\rightarrow] &4&4&4&&&&&\\\none & \none &5&5&5&&&&&\end{ytableau} \\
&&\\
(i) & \begin{ytableau} \ &&&&&1&1&\none&1\\&&&&1&&2&\none&2\\&1&1&1&&2&&\none & &\none[\leftarrow] & \none[\framebox{$4$}_{\,3}] \\1&2&2&2&2&3&3&\none&3\\
&&&&&4&4&\none&4\end{ytableau} &
\begin{ytableau}\none & \none &1 &\none & 1& 1&2&3&3&3&4\\
\none & \none &2&\none &2&3&4&4&4&4&\\
 \none & \none & 4&\none &4&4&&&&&\\
 \none[\framebox{$3$}_{\, 4}] & \none[\rightarrow] & 5&\none &5&5&&&&&\end{ytableau}
&
\begin{ytableau} \none & \none &1 & \none&1& 1&2&3&3&3&4\\
\none & \none &2&\none&2&3&4&4&4&4&\\
\none[3] & \none[\rightarrow] &4&\none&4&4&&&&&\\
\none & \none &5&\none&5&5&&&&&\end{ytableau} \\
&&\\

 (ii)  & \begin{ytableau} \ &&&&&1&1&\none &\none&1\\
 &&&&1&&2&\none &\none&2\\
 &1&1&1&&2&&\none& \none &4 \\
 1&2&2&2&2&3&3&\none &3&\\
&&&&&4&4&\none&4&\end{ytableau} &
\begin{ytableau}\none &\none &1&\none &\none & 1& 1&2&3&3&3&4\\
\none &\none &2&\none &\none&2&3&4&4&4&4&\\
\none &\none & &4&\none&4&4&&&&&\\
\none &\none &3&5&\none &5&5&&&&&\end{ytableau}
&
\begin{ytableau} \none & \none &1 & \none&1& 1&2&3&3&3&4\\
\none & \none &2&\none&2&3&4&4&4&4&\\
\none[3] & \none[\rightarrow] &4&\none&4&4&&&&&\\
\none & \none &5&\none&5&5&&&&&\end{ytableau} \\
&&\\
 (iii) & \begin{ytableau} \ &&&&&1&1&\none &\none&1\\
 &&&&1&&2&\none &\none&2\\
 &1&1&1&&2&&\none& \none &3 \\
 1&2&2&2&2&3&3&\none[\leftarrow] &3&\\
&&&&&4&4&\none&\none &4\end{ytableau} &
\begin{ytableau}\none &\none &1&\none &\none & 1& 1&2&3&3&3&4\\
\none &\none &2&\none &\none&2&3&4&4&4&4&\\
\none &\none & 3&4&\none[\rightarrow]&4&4&&&&&\\
\none &\none &5&\none&\none &5&5&&&&&\end{ytableau}
&
\begin{ytableau}\none &\none &1&\none &\none & 1& 1&2&3&3&3&4\\
\none &\none &2&\none &\none&2&3&4&4&4&4&\\
\none &\none & 3&\none[4]&\none[\rightarrow]&4&4&&&&&\\
\none &\none &5&\none&\none &5&5&&&&&\end{ytableau}
 \\
\end{array} \]

 \ytableausetup{smalltableaux}

\begin{cor} Suppose $T$ is a ptableau, and $\RSK(T) = (PT,T_{\mmax})$. Also assume that $T = \Pf{\tau \choose \omega}$ and $(P,Q) =\RSK{\tau \choose \omega}$ is the classic RSK pair obtained by column insertion of the word $\omega$, starting on the left, using the word $\tau$ to determine the growth tableau $Q$. Then:
\begin{enumerate}
\item $PT = \widehat{P}$. That is, regarding the SSYT $P$ as a (necessarily highest weight) ptableau, its dual $\widehat{P}$ equals the ptableau $PT$ obtained by the ptableaux insertion algorithm $\RSK(T) = (PT,T_{\mmax})$.
\item $Q= T_{\mmax}$. That is, the classic RSK growth tableau $Q$, viewed as a ptableau, is the highest weight element in the irreducible crystal containing $T$.
\item In particular, if $T_{\mmax} = Q$, and $\bw(Q)={\tau \choose \eta}$, then $\eta$ is the highest weight of the word $\omega$ where $\bw(T) = {\tau \choose \omega}$ (noting that since $T$ and $T_{\mmax}$ are in the same crystal, they have the same content and hence the same top row $\tau$ in their biword). In particular, under the \emph{classic} RSK algorithm we can read off the highest weight of $\omega$ from the rows in which the $1$s in $Q$ appear, then the rows in which the $2$s in $Q$ appear, and so on.
\end{enumerate} \label{PT to P corollary}
\end{cor}
\begin{proof}
Statement (1) follows immediately from Lemma~\ref{insert RSK} since by the end of the insertion algorithm we have $P=P^{(k)} = \widehat{PT^{(k)}} = \widehat{PT}$.

For the proof of Statement (2), it is clear that $Q^{(1)} = T_{\mmax}^{(1)}$, and inductively we assume $Q^{(k)} = T_{\mmax}^{(k)}$. Suppose we are computing the ptableau RSK insertion $PT^{(k)} \otimes \framebox{$\tau_{k+1}$}_{\, \omega_{k+1}}$, and that the last violation entry, after it is resolved, becomes $\framebox{$s$}_{\, t}$. Then, since $P^{(k+1)}= \widehat{PT^{(k+1)}}$, and by Statement (1) above, the last box added by RSK insertion on $P$ was $\framebox{$k$}_{\, s}$. In particular, the value $s$ determines the \emph{row} in which the growth in $P^{(k)}$ occurs. The algorithm then augments $T_{\mmax}$ by adding $\framebox{$\tau_{k+1}$}_{\, s}$. This is the same growth we would obtain in $Q^{(k+1)}$ by the classic RSK insertion using the content $\tau_{k+1}$ appearing at the end of the row $s$ in which $P^{(k+1)}$ was augmented.

Lastly, to prove (3), since by (2) the classic RSK $Q$ is actually the highest weight of the crystal containing $\Bw(T)={\tau \choose \omega}$, the \emph{rows} in which the content of $Q$ appears, in the order determined by $\tau$, must be the highest weight of the word $\omega$.
 \end{proof}

 So, in our example,
\[ T = \begin{ytableau} \ & & 1 & 3&4 \\ 1 & 2 & 2& & \\ 3 & 3& 4& 4& \end{ytableau}, \ \
\hbox{ and we have } \Bw(T) = {1122333444  \choose2122331331}. \]
We obtained above that
\[ \RSK(T) = \left(\,  \begin{ytableau} \ & & 1&1&1\\1&1&&&2\\&2&2&2&3\end{ytableau}\ , \ \begin{ytableau} 1&1&2&3&4\\2&3&4&4&\\3&&&&\end{ytableau} \, \right). \]

Thus, were we to have computed the classic RSK on the biword $\Bw(T)$, we could have used the growth tableau

\[ Q = \begin{ytableau} 1&1&2&3&4\\2&3&4&4&\\3&&&&\end{ytableau}  \]
to read off the highest weight of the crystal containing the word
$\omega = (2122331331)$. Noting the rows (left to right) in which the 1's appear in $Q$, then the 2's, etc., we see that the highest weight of the crystal containing $\omega$ must be:
\[ \eta= (1121321221). \]

To emphasize, we see that ptableau insertion $RSK(T) = (PT,T_{\mmax})$ computes $T_{\mmax}$, the highest weight element in the irreducible crystal containing $T$. Computing the ptableau $PT$ not only produces a ptableau in the distinguished crystal that is plactically equivalent to $T$, but the form of $PT$ allows us to compute the \emph{path} of crystal operators from $T$ to $T_{\mmax}$ in this crystal.

\begin{thm} Suppose $T \in \PTab_n$ and $\RSK(T) = (PT,T_{\mmax})$. Let $\alpha_{ij}$ be the number of $i$'s appearing in row $j$ of the ptableau $PT$, and let
\begin{equation} \beta_{st} = \alpha_{1,t} + \cdots + \alpha_{s,t},\end{equation}
and
\begin{equation} e_{(\ell)} = e_{n-\ell}^{\beta_{n-\ell,n}} \cdots e_{1}^{\beta_{1,(1+\ell)}},\end{equation}
and finally
\begin{equation} e_{(*)} = e_{(n-1)} e_{(n-2)} \cdots e_{(1)}. \label{e up to max}\end{equation}
 Then
\begin{equation} T_{\mmax} = e_{(*)}T. \label{e up to lus}\end{equation}
\label{e up formula}
\end{thm}

\begin{proof}
This formula is harder to express notationally than it is to understand. The symbol $\alpha_{ij}$ denotes the number of $i$'s in row $j$ of $PT$, computed from ptableau RSK insertion of a ptableau $T$. By Corollary~\ref{word equiv}, the integers $i$ in row $j$ form a right-contiguous strip in the right-justified $PT^*$, and so we denote the portion of the ptableau $PT$ with this content as $\framebox{$i^{\alpha_{ij}}$}$ to suggest a right-justified string of $i$'s of length $\alpha_{ij}$:

\[\framebox{$i^{\, \alpha_{ij}}$} = \begin{ytableau} \ &\scalebox{0.7}[0.8]{$\cdots$} & \end{ytableau} \hspace{-0.03in}\overbrace{\begin{ytableau}  i & i &\scalebox{0.7}[0.8]{$\cdots$}&i \end{ytableau}}^{\alpha_{ij}}, \]
where the first $\cdots$ denotes a (possibly empty) sequence of blanks in the right-justified form, and the second $\cdots$ denotes a (possibly empty) sequence if $i$'s.

Then, since $PT$ is in the distinguished crystal, we can depict a 5-rowed example in right-justified form as:

\ytableausetup{nosmalltableaux}
\[ PT^* = \begin{ytableau} \ &&&& \scalebox{0.7}[0.8]{$1^{\alpha_{11}}$} \\
&&&\scalebox{0.7}[0.8]{$1^{\alpha_{12}}$} & \scalebox{0.7}[0.8]{$2^{\alpha_{22}}$} \\
& & \scalebox{0.7}[0.8]{$1^{\alpha_{13}}$} & \scalebox{0.7}[0.8]{$2^{\alpha_{23}}$} & \scalebox{0.7}[0.8]{$3^{\alpha_{33}}$} \\
 &  \scalebox{0.7}[0.8]{$1^{\alpha_{14}}$}& \scalebox{0.7}[0.8]{$2^{\alpha_{24}}$}  & \scalebox{0.7}[0.8]{$3^{\alpha_{34}}$} &  \scalebox{0.7}[0.8]{$4^{\alpha_{44}}$} \\
\scalebox{0.7}[0.8]{$1^{\alpha_{15}}$} & \scalebox{0.7}[0.8]{$2^{\alpha_{25}}$} &  \scalebox{0.7}[0.8]{$3^{\alpha_{35}}$} & \scalebox{0.7}[0.8]{$4^{\alpha_{45}}$} & \scalebox{0.7}[0.8]{$5^{\alpha_{55}}$} \end{ytableau}. \]

In our notation, we have

\begin{align*} e_{(1)} &= e_{n-1}^{\beta_{n-1,n}} \cdots e_{3}^{\beta_{3,4}} e_{2}^{\beta_{2}^{2,3}}e_{1}^{\beta_{1,2}} \\
&= e_{n-1}^{\beta_{n-1,n}} \cdots e_{3}^{(\alpha_{14}+\alpha_{24}+\alpha_{34})} e_{2}^{(\alpha_{13}+\alpha_{23})}e_{1}^{\alpha_{12}}. \end{align*}

Note that in rows 1 and 2, we have (but now putting the blocks in \emph{left-justified form}):

\[ ^*PT[1,2] = \begin{ytableau}\ & \scalebox{0.7}[0.8]{$1^{\alpha_{11}}$} \\
\scalebox{0.7}[0.8]{$1^{\alpha_{12}}$} & \scalebox{0.7}[0.8]{$2^{\alpha_{22}}$} \end{ytableau}. \]
Since $\alpha_{11} \geq \alpha_{22}$ there are no $2$'s in row 2 without $1$'s above them. So, the $1$'s in row 2 constitute the right-most content under blanks, we have
\[ e_{1}^{\alpha_{12}} \left(^*PT[1,2]\right) = \begin{ytableau}\scalebox{0.7}[0.8]{$1^{\alpha_{12}}$} & \scalebox{0.7}[0.8]{$1^{\alpha_{11}}$} \\
 & \scalebox{0.7}[0.8]{$2^{\alpha_{22}}$} \end{ytableau}.  \]
But then, in rows 2 and 3 we have:
\[ ^*\left(e_{1}^{\alpha_{12}} PT\right)[2,3] = \begin{ytableau}\
& & \scalebox{0.7}[0.8]{$2^{\alpha_{22}}$} \\
  \scalebox{0.7}[0.8]{$1^{\alpha_{13}}$} & \scalebox{0.7}[0.8]{$2^{\alpha_{23}}$} & \scalebox{0.7}[0.8]{$3^{\alpha_{33}}$} \end{ytableau},  \]
and again, the word condition $\alpha_{22} \geq \alpha_{33}$ implies the $3$'s in row $3$ must be covered by the $2$'s in row $2$. Thus, the 1's and 2's in row 3 constitute the right-most content under blanks in those rows. But then
\[ ^*e_{2}^{\beta_{2,3}}e_{1}^{\beta_{1,2}}\left(PT[2,3]\right)=
{^*e_{2}^{(\alpha_{13}+\alpha_{23})}\left(e_{1}^{\alpha_{12}} PT\right)[2,3]} = \begin{ytableau}\scalebox{0.7}[0.8]{$1^{\alpha_{13}}$}
&\scalebox{0.7}[0.8]{$2^{\alpha_{23}}$} & \scalebox{0.7}[0.8]{$2^{\alpha_{22}}$} \\
   &  & \scalebox{0.7}[0.8]{$3^{\alpha_{33}}$} \end{ytableau} \]
Thus, we claim, the overall effect of applying $e_{(1)}$, then $e_{(2)}$, etc., to $PT$ becomes:

\[ PT^* = \begin{ytableau} \ &&&& \scalebox{0.7}[0.8]{$1^{\alpha_{11}}$} \\
&&&\scalebox{0.7}[0.8]{$1^{\alpha_{12}}$} & \scalebox{0.7}[0.8]{$2^{\alpha_{22}}$} \\
& & \scalebox{0.7}[0.8]{$1^{\alpha_{13}}$} & \scalebox{0.7}[0.8]{$2^{\alpha_{23}}$} & \scalebox{0.7}[0.8]{$3^{\alpha_{33}}$} \\
 &  \scalebox{0.7}[0.8]{$1^{\alpha_{14}}$}& \scalebox{0.7}[0.8]{$2^{\alpha_{24}}$}  & \scalebox{0.7}[0.8]{$3^{\alpha_{34}}$} &  \scalebox{0.7}[0.8]{$4^{\alpha_{44}}$} \\
\scalebox{0.7}[0.8]{$1^{\alpha_{15}}$} & \scalebox{0.7}[0.8]{$2^{\alpha_{25}}$} &  \scalebox{0.7}[0.8]{$3^{\alpha_{35}}$} & \scalebox{0.7}[0.8]{$4^{\alpha_{45}}$} & \scalebox{0.7}[0.8]{$5^{\alpha_{55}}$} \end{ytableau}
\stackrel{e_{(1)}}{\longrightarrow}
\begin{ytableau} \ &&&\scalebox{0.7}[0.8]{$1^{\alpha_{12}}$}& \scalebox{0.7}[0.8]{$1^{\alpha_{11}}$} \\
&&\scalebox{0.7}[0.8]{$1^{\alpha_{13}}$} & \scalebox{0.7}[0.8]{$2^{\alpha_{23}}$} & \scalebox{0.7}[0.8]{$2^{\alpha_{22}}$} \\
&  \scalebox{0.7}[0.8]{$1^{\alpha_{14}}$}& \scalebox{0.7}[0.8]{$2^{\alpha_{24}}$}  & \scalebox{0.7}[0.8]{$3^{\alpha_{34}}$} & \scalebox{0.7}[0.8]{$3^{\alpha_{33}}$} \\
\scalebox{0.7}[0.8]{$1^{\alpha_{15}}$} & \scalebox{0.7}[0.8]{$2^{\alpha_{25}}$} &  \scalebox{0.7}[0.8]{$3^{\alpha_{35}}$} & \scalebox{0.7}[0.8]{$4^{\alpha_{45}}$}  &  \scalebox{0.7}[0.8]{$4^{\alpha_{44}}$} \\
&&&& \scalebox{0.7}[0.8]{$5^{\alpha_{55}}$} \end{ytableau}
\stackrel{e_{(2)}}{\longrightarrow}
\begin{ytableau} \ &&\scalebox{0.7}[0.8]{$1^{\alpha_{13}}$}&\scalebox{0.7}[0.8]{$1^{\alpha_{12}}$}& \scalebox{0.7}[0.8]{$1^{\alpha_{11}}$} \\
&\scalebox{0.7}[0.8]{$1^{\alpha_{14}}$}& \scalebox{0.7}[0.8]{$2^{\alpha_{24}}$} & \scalebox{0.7}[0.8]{$2^{\alpha_{23}}$} & \scalebox{0.7}[0.8]{$2^{\alpha_{22}}$} \\
\scalebox{0.7}[0.8]{$1^{\alpha_{15}}$} & \scalebox{0.7}[0.8]{$2^{\alpha_{25}}$} &  \scalebox{0.7}[0.8]{$3^{\alpha_{35}}$}& \scalebox{0.7}[0.8]{$3^{\alpha_{34}}$} & \scalebox{0.7}[0.8]{$3^{\alpha_{33}}$} \\
&& & \scalebox{0.7}[0.8]{$4^{\alpha_{45}}$}  &  \scalebox{0.7}[0.8]{$4^{\alpha_{44}}$} \\
&&&& \scalebox{0.7}[0.8]{$5^{\alpha_{55}}$} \end{ytableau}
\stackrel{\cdots}{\longrightarrow}
\begin{ytableau} \ \scalebox{0.7}[0.8]{$1^{\alpha_{15}}$} &\scalebox{0.7}[0.8]{$1^{\alpha_{14}}$}&\scalebox{0.7}[0.8]{$1^{\alpha_{13}}$}&\scalebox{0.7}[0.8]{$1^{\alpha_{12}}$}& \scalebox{0.7}[0.8]{$1^{\alpha_{11}}$} \\
&\scalebox{0.7}[0.8]{$2^{\alpha_{25}}$}& \scalebox{0.7}[0.8]{$2^{\alpha_{24}}$} & \scalebox{0.7}[0.8]{$2^{\alpha_{23}}$} & \scalebox{0.7}[0.8]{$2^{\alpha_{22}}$} \\
 &  &  \scalebox{0.7}[0.8]{$3^{\alpha_{35}}$}& \scalebox{0.7}[0.8]{$3^{\alpha_{34}}$} & \scalebox{0.7}[0.8]{$3^{\alpha_{33}}$} \\
&& & \scalebox{0.7}[0.8]{$4^{\alpha_{45}}$}  &  \scalebox{0.7}[0.8]{$4^{\alpha_{44}}$} \\
&&&& \scalebox{0.7}[0.8]{$5^{\alpha_{55}}$} \end{ytableau}. \]
This final form has all $i$'s in row $i$, and hence is the highest weight element in the distinguished crystal containing $PT$.

We only need argue that, at each stage, the appropriate sequence of $e_i$ in a given $e_{(k)}$ is non-null. Suppose we have applied $e_{(\ell-1)} \cdots e_{(2)} e_{(1)}$ to $PT$. The next step would be to apply $e_{(\ell)} = e_{n-\ell}^{\beta_{n-\ell,n}} \cdots e_{1}^{\beta_{1,(1+\ell)}}$. We suppose, inductively, that we have applied the initial factors: $e_{k}^{\beta_{k,k+\ell}}\cdots e_{2}^{\beta_{2,2+\ell}}e_{1}^{\beta_{1,1+\ell}}$, assuming that rows $k-1$, $k$, and $k+1$ have the form:

\[ \begin{array}{|c|c|c|c|c||c|c|c|} \hline
 & 1^{\alpha_{1,\ell +k}} & 2^{\alpha_{2,\ell+k}} & \dots & k^{\alpha_{k,\ell+k}} & k^{\alpha_{k, k+\ell -1}} & \dots & k^{\alpha_{k,k}} \\ \hline
  & &&&&(k+1)^{\alpha_{k+1,\ell +k}} & \dots & (k+1)^{\alpha_{k+1,k+1}} \\ \hline
  1^{\alpha_{1,\ell+k+1}} & 2^{\alpha_{2,\ell+k+1}}&3^{\alpha_{3,\ell+k+1}}& \dots &
  (k+1)^{\alpha_{(k+1),\ell+k+1}} & (k+2)^{\alpha_{(k+2),\ell+k+1}} & \dots & (k+2)^{\alpha_{(k+2),k+2}} \\ \hline
  \end{array} \]

(The double bars are just for clarity, to distinguish the portion of the ptableau $PT$ that is not to be moved on the right, and the portions that will be moved by subsequent operators on the left.) The next step is to apply
\[ e_{k+1}^{\beta_{k+1,k+1+\ell}} = e_{k+1}^{ \alpha_{1,\ell+k+1} + \cdots + \alpha_{(k+1),\ell+k+1}}. \]

Since $PT$ satisfies the word condition, we know that the number of $(k+1)$'s in rows $k+1$ through $k+\ell$ are greater than or equal to the number of $(k+2)$'s in rows $k+2$ through $k+\ell+1$. Thus, in rows $k$ and $k+1$ above, to the right of the double line, every $k+2$ is covered by some $k+1$ above it, so that all the content in row $k+1$ with a blank above it lies to the \emph{left} of the double line, resulting in all of the content to the left of the line being moved up into row $k$ by the application of $e_{k+1}^{\beta_{k+1,k+1+\ell}}$.

Thus, the formulas given above take $PT$ to its associated highest weight.
Since $T$ and $PT$ are plactically equivalent by Corollary~\ref{PT to P corollary}, this same sequence of crystal operators takes $T$ to $T_{\mmax}$.
\end{proof}

Let's consider an example. Suppose \ytableausetup{smalltableaux}
\[T = \ytableaushort{\none\none11\none\none,\none\none2344,11345\none}* {6,6,6}= \Pf{\tau \choose \omega} =\Pf{11112334445 \choose 33112323223}. \]
We calculate:
\[\RSK(T)=(PT,T_{\mmax}) = \left( \,\begin{ytableau} \ & & &&1&1\\ &&1&1&2&2\\1&1&&2&3&3 \end{ytableau}\,, \, \begin{ytableau} 1 & 1&1&1&4&4 \\ 2 &3&5&&& \\ 3&4&&&& \end{ytableau}\,\right)\]
and
\[ \RSK{\tau \choose \omega} = (P,Q) = \left( \, \begin{ytableau} 1&1&2&2&3&3\\2&2&3&&&\\3&3&&&& \end{ytableau}, \, \begin{ytableau}1 & 1&1&1&4&4 \\ 2 &3&5&&& \\ 3&4&&&& \end{ytableau}
\, \right) . \]

We compute a crystal operator path from $PT$ to its highest weight, moving $1$'s out of row $2$, then the $1$'s and $2$'s out of row $3$, and then the remaining $1$'s out of row $2$:
\begin{align*} e_{(*)}PT = e_1^{2} e_{2}^{3} e_{1}^2\ &\begin{ytableau} \ & & &&1&1\\ &&1&1&2&2\\1&1&&2&3&3 \end{ytableau}\\
 =e_1^{2} e_{2}^{3}\ & \begin{ytableau} \ & &1 &1&1&1\\ &&&&2&2\\1&1&&2&3&3 \end{ytableau} \\
 = e_1^{2}\  & \begin{ytableau} \ & &1 &1&1&1\\ 1&1&&2&2&2\\&&&&3&3 \end{ytableau}= \begin{ytableau} 1 & 1&1 &1&1&1\\ 2&2&2&&&\\3&3&&&& \end{ytableau}
\end{align*}

But then, in accordance with Theorem~\ref{e up formula}:
\vspace{-0.1in}
\begin{align*} e_{(*)}T = e_1^{2} e_{2}^{3} e_{1}^2\ & \ytableaushort{\none\none11\none\none,\none\none2344,11345\none}* {6,6,6} \\
=e_1^{2} e_{2}^{3}\ & \ytableaushort{\none\none1144,\none\none23\none\none,11345}* {6,6,6}\\
=e_1^{2}\  &\ytableaushort{\none\none1144,11235,34}* {6,6,6}
=\ \ytableaushort{111144,235,34}* {6,6,6}= T_{\mmax}, \end{align*}
which is highest weight.

(We note that the notation $\alpha_{ij}$ in the statement of Theorem~\ref{e up formula} is also the same as the $(i,j)$ entry of the \emph{matrix} $M$ associated to the ptableau $T$ (or its biword $\bw{T} = {\tau \choose \omega}$). Thus, for matrices in the distinguished crystal the above formula also computes the path to highest weight, though we do not pursue this here.)

\medskip

\subsection{Summary}

Our results in Corollary~\ref{PT to P corollary} (2) and Theorem~\ref{e up formula} demonstrate that the ptableau approach to RSK, for the purposes of analyzing crystal graphs, is quite natural and possesses advantages over the classical RSK approach. Classically, the RSK algorithm converts a biword $\tau \choose \omega$ into a SSYT pair $(P,Q)$. We adopt the view that (1) we should regard the SSYT $(P,Q)$ in RSK as perforated tableaux, and (2) we should actually replace the pair $(P,Q)$ (and the duality $\widehat{(P,Q)} = (Q,P)$ ) with the ptableau pair $(\widehat{P},Q)$ (and the duality $\widehat{(\widehat{P},Q)} = (\widehat{Q}, P)$ ) as a more natural invariant with which to study plactic equivalence in crystal graphs. Similar to the classic RSK case, we obtain the commutative diagram:

\[ \begin{array}{cccl} (\widehat{P},Q) & \stackrel{e_i}{\rightarrow} & (e_{i}\widehat{P},Q)& \hbox{$\Leftarrow$ In the same crystal} \\ \raisebox{-0.05in}{$\widehat{\ }$} \downarrow &&\raisebox{-0.05in}{$\widehat{\ }$} \downarrow& \\
(Q,\widehat{P}) & \stackrel{e_i}{\rightarrow} & (Q,e_{i}\widehat{P}) & \hbox{$\Leftarrow$ Plactically equivalent} \end{array} \]
(where now $e_{i}\widehat{P}$ is defined as a ptableau crystal operation). However, since $\widehat{P} = PT$ and $Q=T_{\mmax}$ under the ptableaux insertion $\RSK(T) = (PT,T_{\mmax}) = (\widehat{P},Q)$ we immediately conclude:

\begin{thm} If $T$ and $T'$ are ptableaux, then $T \equiv T'$ if and only if $\widehat{T}$ and $\widehat{T'}$ lie in the same irreducible crystal graph.
\end{thm}

And, in fact,
\begin{itemize}
\item  If $RSK(T) = (PT,T_{\mmax})$ and $RSK(T') = (PT', T_{\mmax}')$, then $PT=PT'$ not only implies $T \equiv T'$, but $PT$ is a common ptableau \emph{plactically equivalent} to $T$ and $T'$. Furthermore, by Theorem~\ref{e up formula} we can read off from $PT$ a path from (either) $T$ or $T'$ to their respective highest weights.
\item Certainly $T_{\mmax} = T_{\mmax}'$ implies $T$ and $T'$ lie in the same irreducible crystal, for the now obvious reason that $T_{\mmax}$, computed by the insertion process, \emph{is the common highest weight element} of the irreducible crystal graph in which both are found. 
\end{itemize}

\section{Algorithms for Lusztig Involution on Ptableaux Crystals}

Suppose $T_{\mmax}$ is some highest weight ptableau (hence, a SSYT by Theorem~\ref{highest weight ptab}). In the authors'~\cite{ptab} a combinatorial algorithm for computing the corresponding \emph{lowest weight} element $T_{\min}$ was given in terms of \emph{evacuation}, $\Evac(T_{\mmax})$, a key step in the Sch\"utzenberger involution.

For example, given a SSYT $T_{\mmax}$:
\ytableausetup{smalltableaux}
 \[  \begin{ytableau}
1 & 1 & 2 & 2 & 3 &4 \\
2 & 3 & 3 & 4& & \\
3 &4 & 5 & & & \end{ytableau}, \]
we compute $\Evac(T_{\mmax})$ by first performing inward jeu de taquin on the outer corner in row 2:
 \[  \begin{ytableau}
1 & 1 & 2 & 2 & 3 &4 \\
2 & 3 & 3 & 4& & \\
3 &4 & 5 & & & \end{ytableau} \Rightarrow \begin{ytableau}
1 & 1 & 2 & 2 & 3 &4 \\
2 & 3 & 3 & &4 & \\
3 &4 & 5 & & & \end{ytableau} \Rightarrow \begin{ytableau}
1 & 1 & 2 & 2 & 3 &4 \\
2 & 3 &   & 3  &4  & \\
3 &4 & 5 &     &    & \end{ytableau} \Rightarrow \begin{ytableau}
1 & 1 & 2 & 2 & 3 &4 \\
2 &   & 3  & 3  &4  & \\
3 &4 & 5 &     &    & \end{ytableau} \Rightarrow \begin{ytableau}
1 & 1 & 2 & 2 & 3 &4 \\
 &  2 & 3  & 3  &4  & \\
3 &4 & 5 &     &    & \end{ytableau}
 \Rightarrow \begin{ytableau}
 \ & 1 & 2 & 2 & 3 &4 \\
 1& 2  &3  & 3  &4  & \\
3 &4 & 5 &     &    & \end{ytableau}. \]
Continuing inner jeu de taquin for all outer corners, we obtain:

\[ \Evac(T_{\mmax}) =\begin{ytableau}
 \ &    &  &  1 &  2 & 3 \\
   &   &   2  &  2  &  3  &4 \\
1 &  3& 3 &  4  & 4 & 5 \end{ytableau}. \]


In the authors'~\cite{ptab}, Theorem 8.5, it was proved that:
\[ T_{\min} = \Evac(T_{\mmax}). \]

 We now present a more general result. As an application of our results above, we give \emph{two} different methods for combinatorially computing the Lusztig involution $T_{\Lus}$ of an arbitrary ptableau $T$. One method uses our insertion and an uninsertion algorithm, the second, the formula of Theorem~\ref{e up formula} (also based on ptableau RSK insertion) to compute $T_{\Lus}$.

We first need to describe the appropriate \emph{uninsertion} algorithm on ptableaux.  It will be clear from construction that this process is the inverse of the RSK algorithim in Theorem \ref{RSK algorithm}.

\bigskip

\noindent {\bf\texttt{Algorithm} $\RSK^{-1}$}\\

\noindent Input data: Ptableaux $PT, T_{\mmax} \in \PTab_{(m,n)}^k$, where $PT$ is in the distinguished crystal, and $T_{\mmax}$ is highest weight.

\noindent Output: A biword $\binom{\tau}{\omega}\in \Bw([m]\times [n])^{\otimes k}$, or equivalently, its ptableau $T$, such that $\RSK(T)=(PT,T_{\mmax})$.\\

\noindent Initialization: Set $j=1$, $PT^{<j>}=PT$, $T_{\mmax}^{<j>} = T_\mmax$,  and $\binom{\tau}{\omega}$ to the empty biword.
\begin{enumerate}

\item  Set:
\begin {enumerate}
\item  $\ell_j$ to the number of columns in $PT^{<j>}$,
\item $\framebox{t}_{\, a}$ to the largest, right-most entry in $T_\mmax^{<j>}$, with $a$ the row containing this cell,
\item $T_\mmax^{j+1} = T_\mmax^{<j>}$ with $\framebox{t}_{\, a}$ removed.
\end{enumerate}

\noindent
 Let $\framebox{a}_{\, r}$ be the left-most $a$ in $PT^{<j>}$ (with $r$ its row).  Set $c$ to the number of the column containing $\framebox{a}_{\, r}$ (with the left-most column being column 1).  We recursively define the ``uninsertion" process in $PT^{<j>}$:\\

\item Consider column $c+1$, the column immediately to the right of the column containing $\framebox{a}_{\, r}$.
\begin{enumerate}
\item If this column is empty, meaning $c+1> \ell_j$ then update $\binom{\tau}{\omega}$ to be $\binom{t}{a}\otimes \binom{\tau}{\omega}$, increase $j$ by $1$, and if $j < k$, return to Step 1.
\item If this column is not empty, then by the word condition and right justification, if the column contains some number $i$ (for $i>1$), then it contains $i-1$.  Let $\framebox{b}_{\, s}$ be the lowest cell in column $c+1$ that is in row $r$ or higher.
\begin{enumerate}
\item Remove $\framebox{b}_{\, s}$ from column $c+1$.
\item Move $\framebox{a}_{\, r}$ to column $c+1$ (so one column to its right, but still in row $r$).
\item Replace the content $a$ with $b$.
\item Update $c$ to $c+1$.
\item Return to Step 2,  with $\framebox{a}_{\, r}$ updated to $\framebox{b}_{\, s}$.
\end{enumerate}
\end{enumerate}
\end{enumerate}
This process can be described visually by the diagram below:
\begin{figure}[H]
\begin{center}
{\includegraphics[scale=.27]{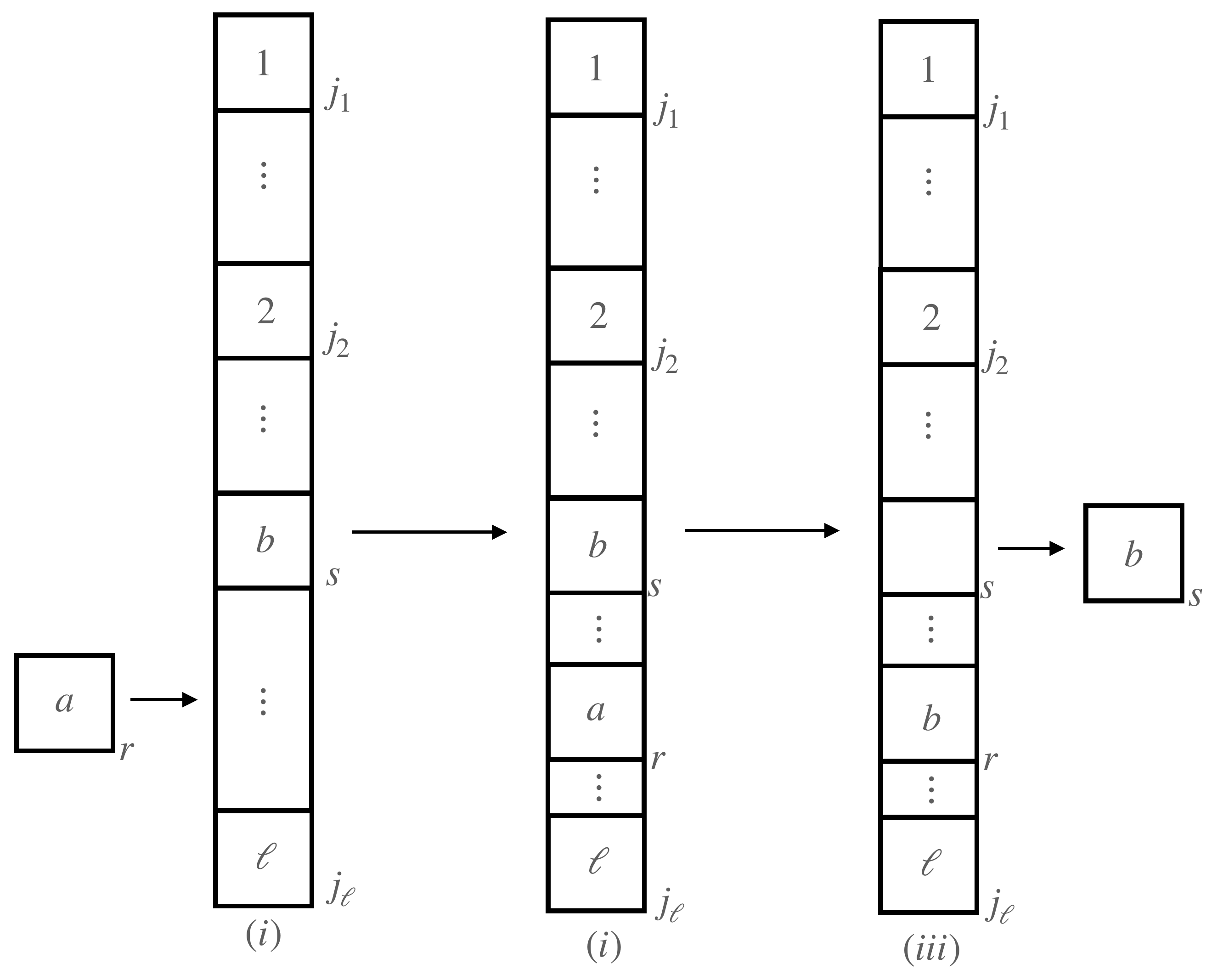}}
\end{center}
\end{figure}
\noindent
with the process then recursively performed, moving $\framebox{b}_{\, s}$ into the column immediately to its right.

By an inspection of the proof of Theorem~\ref{PT to P corollary}, the following is easily concluded:

\begin{thm} Let $PT \in \PTab_{n}$ lie in the distinguished crystal, and let $T_{\mmax} \in \PTab_{n}$ be highest weight, with the same shape at $PT$ as a SSYT. Then the map $\RSK^{-1}:(PT,T_{\mmax})\mapsto T$ given above inverts the ptableaux insertion algorithm. That is, $\RSK^{-1}\circ RKS$ is the identity on ptableaux.
\end{thm}

With this, we now give two different methods for computing Lusztig involutions on ptableaux crystals.

\begin{df} Let $T$ be a node in a fixed irreducible crystal ${\cal B}$ of type $A_{n-1}$ We say a node $T_{\Lus} \in {\cal B}$ is the \emph{Lusztig involution} of $T$ when: If $T_{\mmax}$ is the highest weight element of ${\cal B}$ and $e_{i_1} \cdots e_{i_{\ell}} T = T_{\mmax}$ is a crystal operator path taking $T$ to $T_{\mmax}$, then if $T_{\mmin}$ is the lowest weight element of ${\cal B}$, we have
\[ e_{n-i_1} \cdots e_{n-i_{\ell}} T_{\mmin} = T_{\Lus}. \]
\end{df}

To obtain our first algorithm to compute $T_{\Lus}$ from $T$, we need a notion of duality different from that used above. The form we study here has appeared in many places, including the excellent text by Bump and Schilling~\cite{BumpSchilling}. To avoid confusion with the definition of dual used previously, we adopt the following convention:

\begin{df}[See~\cite{BumpSchilling}, p.\ 15] Given an irreducible type $A_{n-1}$ crystal ${\cal B}$, we say a crystal ${\cal B}'$ is \emph{BS-dual} to ${\cal B}$ if there is a bijection $\phi: {\cal B} \rightarrow {\cal B}'$ such that
\[ f_i \phi(T) = \phi( e_{n-i} T), \quad \hbox{and} \quad e_i \phi(T) = \phi( f_{n-i} T) \]
for all nodes $T \in {\cal B}$.
\end{df}

If ${\cal B}$ is an irreducible crystal of \emph{ptableaux}, it is immediate from Proposition~\ref{rot prop} that the image of ${\cal B}$ under the map $\Rot$ is BS-dual to ${\cal B}$. Set ${\cal B}' = \Rot({\cal B})$. Then ${\cal B}$ and ${\cal B}'$ are isomorphic, irreducible crystals (even though the map $\Rot$ is definitely \emph{not} a crystal isomorphism) and, given $T \in {\cal B}$, the image $\Rot(T) \in {\cal B}'$ is a ptableau that is plactically equivalent to $T_{\Lus} \in {\cal B}$. That is, $T_{\Lus}$ and $\Rot(T)$ live in the same ``location'' in their respective crystals. We exploit this in the following result.

\begin{thm} Let $T \in \PTab_{n}$ and let $RSK(T) = (PT,T_{\mmax})$. Let $T_{\Lus}$ denote the Lusztig involution of $T$. Then
\[ T_{\Lus} = \Rot(\RSK^{-1}(PT,Rot(\Evac(T_{\mmax})))). \]
\label{insert formula}
\end{thm}
\vspace{-0.2in}
\begin{proof} Most of the work has been argued above. From~\cite{ptab} we know that $\Evac(T_{\mmax})=T_{\mmin}$, the lowest weight element in the irreducible crystal containing $T$. Consequently, $\Rot(\Evac(T_{\mmax}))$ is the highest weight element in the BS-dual crystal. Thus, the ptableau $\RSK^{-1}(PT,\Rot(\Evac(T_{\mmax})))$ is plactically equivalent to $T$ (since they share the same $PT$ under the ptableaux RSK correspondence), and hence applying $\Rot$ to this ptableau maps back to ${\cal B}$, so that the image must be $T_{\Lus}$.
\end{proof}

An alternate method to compute $T_{\Lus}$ from a ptableau $T$,
we have the following:

\begin{thm} Suppose $T \in PTab_n$ and $RSK(T) = (PT,T_{\mmax})$. Let $e_{(*)}$ be defined as in Theorem~\ref{e up formula}.

Let $\overbrace{e_{(*)}}$ denote the product of crystal operators formed by writing the factors in the reverse order as they appear in $e_{(*)}$, and replacing each $e_{i}$ with $e_{n-i}$. Then we have:
\[ \overbrace{e_{(*)}}\, \Evac(T_{\mmax}) = \overbrace{e_{(*)}}\, T_{\mmin}=T_{\Lus}. \]
\end{thm}

Using the example following Theorem~\ref{e up formula}, if

\ytableausetup{smalltableaux}
\[T = \ytableaushort{\none\none11\none\none,\none\none2344,11345\none}* {6,6,6} \ \ \hbox{then} \ \
\RSK(T)=(PT,T_{\mmax}) = \left( \,\begin{ytableau} \ & & &&1&1\\ &&1&1&2&2\\1&1&&2&3&3 \end{ytableau}\,, \, \begin{ytableau} 1 & 1&1&1&4&4 \\ 2 &3&5&&& \\ 3&4&&&& \end{ytableau}\,\right). \]
By Theorem~\ref{e up formula} above we have $e_{(*)}=e_1^{2} e_{2}^{3} e_{1}^2$ and hence we compute $T_{\Lus}$ from $T$ by:

\[ \overbrace{e_1^{2} e_{2}^{3} e_{1}^2} \, \Evac(T_{\mmax}) = e_2^{2} e_{1}^{3} e_{1}^2\, \begin{ytableau} \ & & &&1&1\\ &&&2&3&4\\1&1&3&4&4&5 \end{ytableau} = \begin{ytableau} \ &1 &1 &1&4&4\\ 1&&2&3&&5\\&&3&4&& \end{ytableau} =T_{\Lus}\, . \]

We summarize these two approaches pictorially below:

\[ \ytableausetup{smalltableaux}
\begin{tikzpicture}
\node (A) { \ytableaushort{111144,235\none\none\none,34}* {6,6,6}$\,=T_{\mmax}$};
\node [above of=A,node distance = 2cm](G) {\hspace{-0.5in}$\cal B$};
\node [right of=A,node distance = 9cm](B) {\ytableaushort{122355,234,55}* {6,6,6}$\,=T_{\mmax}'$};
\node [above of=B,node distance = 2cm](G) {\hspace{-0.4in}$\Rot({\cal B})=\cal B'$};
\node [left of=A,below of =A,node distance = 3cm](C) {\ytableaushort{\none\none\none\none11,\none\none2344,1134\none5}* {6,6,6}$\,=T$};
\node [below of =C,node distance = 3cm](D) {\ \ \ \  \ytableaushort{\none11144,1\none23\none5,\none\none34}* {6,6,6}$\,=T_{\Lus}$};
\node [below of =A,node distance = 9cm](E) {\ytableaushort{\none\none\none\none11,\none\none\none234,113445}* {6,6,6}$=\,T_{\mmin}$};
\node [below of =B,left of=B,node distance = 3cm](F) {\hspace{1.3in} \ytableaushort{\none\none\none\none23,\none\none1345,22555}* {6,6,6}$\,= \RSK^{-1}(PT,T_{\mmax}')$};
 \node[ellipse,
    draw = black,
    minimum width = 8.5cm,
    minimum height = 11.2cm] (e) at (-0.5,-4.5) {};
  \node[ellipse,
    draw = black,
    minimum width = 8cm,
    minimum height = 11.2cm] (e) at (8.5,-4.5) {};
 \draw[->] (C.90) to node [above] {$e_{(*)}\ \ \ $} (A.200);
  \draw[->] (F.90) to node [above] {$e_{(*)}\ \ \ $} (B.200);
  \draw[->] (A.270) to [out=-20,in=70]node [right] {$\Evac$} (E.130);
  \draw[->] (E.-45) to [out=-20,in=-20]node [right] {$\ \Rot$} (B.-60);
   \draw[->] (C.-45) to [out=-20,in=40]node [right] {$\ \Lus$} (D.80);
    \draw[->] (F.-150) to [out=-140,in=0]node [below] {$\Rot$} (D.0);
    \draw[->] (E.165) to node [below] {$\overbrace{e_{(*)}}\ \ \ $} (D.270);
   \node [above of =A,node distance = 3cm](G) {$\RSK(T)=(PT,T_{\mmax}) = \left( \,\begin{ytableau} \ & & &&1&1\\ &&1&1&2&2\\1&1&&2&3&3 \end{ytableau}\,, \, \begin{ytableau} 1 & 1&1&1&4&4 \\ 2 &3&5&&& \\ 3&4&&&& \end{ytableau}\,\right)$};
 \draw[->] (C.135) to [out=120,in=-150]node [left] {$\ \RSK$} (G.200);
\end{tikzpicture}
 \]

\noindent \mbox{} \hrulefill \mbox{}


\begin{thebibliography}{99}

\bibitem{ptab} G.\ Appleby and Tamsen Whitehead:  Perforated Tableaux: A Combinatorial Model for Crystal Graphs in Type $A_{n-1}$. Algebr Represent Theor (2022). https://doi.org/10.1007/s10468-022-10135-4



\bibitem{BumpSchilling} D.\ Bump and A.\ Schilling: Crystal Bases, Representations and Combinatorics. World Scientific, Hackensack, NJ (2017)






\bibitem{Fulton} W.\ Fulton. Young Tableaux. Cambridge University Press (1997)

\bibitem{GerberLecouvey} T.\ Gerber and C.\ Lecouvey. Duality and bicrystals on infinite binary matrices. ArXiv: 2009.10397 (2020)





\bibitem{Kash-Nak} M.\ Kashiwara and T.\ Nakashima: Crystal Graphs for representations of the $q$-analogue of classical Lie algebras, \emph{J. Algebra}, {\bf 165}, 2, pp. 295-345 (1994)

and semistandard oscillating tableaux. ArXiv:1910.04459 (2019)




\bibitem{Shimozono} M.\ Shimozono. Crystals for Dummies. https://www.aimath.org/WWN/kostka/crysdumb.pdf (2005)

\bibitem{vanLeeuwen} M.\ van Leeuwen. Double Crystals of Binary and Integral Matrices. \emph{Electronic Jounal of Combinatorics}. 13 (2006)


\end{thebibliography}
\end{document}